\documentclass[review]{elsarticle}
\usepackage{graphicx}
\usepackage{epstopdf}
\usepackage{lineno,hyperref}
 \usepackage{geometry}
\usepackage{epsfig}
\usepackage{lineno,hyperref}
\usepackage{bookmark}
\usepackage{multirow}
 \usepackage{amsfonts}
\usepackage{multirow}
 \usepackage{algorithm}
 \usepackage[figuresright]{rotating}
 \usepackage{amscd}
\usepackage{mathrsfs}
\usepackage{booktabs,longtable}
\usepackage{indentfirst}
\usepackage{subfigure}
\usepackage{amstext}
\usepackage{amssymb}
\usepackage{fancyhdr}
\usepackage{amsmath}
\usepackage{color}
\usepackage{bm}
\usepackage{algorithmicx}
\usepackage{algpseudocode}
\usepackage{amsmath}
\modulolinenumbers[5]
\biboptions{sort&compress}
\newtheorem{definition}{Definition}[section]

\newtheorem{theorem}{Theorem}[section]
\newtheorem{corollary}{Corollary}[section]
\newtheorem{lemma}{Lemma}[section]
\newtheorem{proposition}{Proposition}[section]
\newtheorem{remark}{Remark}[section]
\newtheorem{example}{Example}[section]

\makeatletter
\@addtoreset{equation}{section}
\makeatother

\journal{Journal of \LaTeX\ Templates}

\begin{document}

\begin{frontmatter}

\title{Sketch-and-project methods for tensor linear systems}

\tnotetext[mytitlenote]{The work is supported by the National Natural Science Foundation of China (No. 11671060) and the Natural Science Foundation Project of CQ CSTC (No. cstc2019jcyj-msxmX0267)}

\author{Ling Tang, Yajie Yu, Yanjun Zhang, \ Hanyu Li\corref{mycor}}
\cortext[mycor]{Corresponding author. E-mail addresses: lihy.hy@gmail.com or hyli@cqu.edu.cn.}

\address{College of Mathematics and Statistics, Chongqing University, Chongqing 401331, P.R. China}
\begin{abstract} For tensor linear systems with respect to the popular t-product, we first present the sketch-and-project method and its adaptive variants. Their Fourier domain versions are also investigated. Then, considering that the existing sketching tensor or way for sampling has some limitations, we propose two improved strategies. Convergence analyses for the methods mentioned above are provided. We compare our methods with the existing ones using synthetic and real data. 
 Numerical results show that they have quite decent performance  in terms of the number of iterations and running time.
\end{abstract}

\end{frontmatter}


\section{Introduction} \label{secint}
In this paper, we aim to solve the following consistent tensor linear systems
\begin{equation}\label{sec1e1}
\mathcal{A}\ast\mathcal{X}=\mathcal{B},
\end{equation}
where $\mathcal{A}\in \mathbb{R}^{m\times n\times l}$ , $\mathcal{X}\in \mathbb{R}^{n\times p\times l}$  and $\mathcal{B}\in \mathbb{R}^{m\times p\times l}$ are third-order tensors, and the operator $*$ denotes the t-product introduced by Kilmer and Martin \cite{kilmer2011factorization}. The problem (\ref{sec1e1}) arises in many applications including tensor dictionary learning \cite{ soltani2016tensor, newman2020nonnegative}, tensor neural network \cite{newman2018stable}, boundary finite element method \cite{ahn1991efficient, alavikia2011electromagnetic, czuprynski2012parallel}, etc. 
For t-product, it has an advantage that it can reserve the information inherent in the flattening of a tensor and, 
with it, many properties of numerical linear algebra can be extended to third and high order tensors \cite{braman2010third, jin2017generalized, lund2020tensor, miao2020generalized, miao2021t,  zheng2021t, qi2021t}. Hence, extensive works on t-product have appeared in recent years and have also been applied in many areas such as image and signal processing \cite{kilmer2013third, soltani2016tensor, tarzanagh2018fast}, computer vision \cite{xie2018unifying, yin2018multiview}, data denoising \cite{zhang2018nonlocal}, low-rank tensor completion\cite{semerci2014tensor,  zhang2014novel, zhang2016exact, zhou2017tensor}, etc. We will review the basic knowledge on t-product in Section \ref{sec:pre}. 

 To solve the problem (\ref{sec1e1}), Ma and Molitor \cite{ma2021randomized} extended the matrix randomized Kaczmarz (MRK) method \cite{kaczmarz1, Strohmer2009} and called it the tensor randomized Kaczmarz (TRK) method. Later, this method was applied to tensor recovery problem \cite{chen2021regularized}. Recently, Du and Sun extended the matrix randomized extended Kaczmarz 
method to the inconsistent tensor recovery problem \cite{du2021randomized}. As we know, the MRK method is a popular iterative method for solving large-scale matrix linear systems, i.e., the case for $l=1$ and $p=1$ in the problem (\ref{sec1e1}), and it has wide developments; see for example \cite{needell2014paved, liu2016accelerated, nutini2016convergence, de2017sampling, Bai2018, Bai2018r, haddock2021greed}.
Most of these methods can be unified into 
 the sketch-and-project (MSP) method and its adaptive variants proposed by Gower et al. \cite{gower2015randomized,gower2019adaptive}.
 Inspired by the above research, we propose the tensor sketch-and-project (TSP) method and its adaptive variants to solve the problem (\ref{sec1e1}), followed by their theoretical guarantees. Meanwhile, we also present their Fourier domain versions and analyze the corresponding convergence. So,
 the TRK method and its theoretical analysis \cite{ma2021randomized} will be the special case of our results.

Besides the randomized algorithms in \cite{ma2021randomized,chen2021regularized,du2021randomized} mentioned above, there are some research based on random sketching technique for t-product; see for example \cite{tarzanagh2018fast,  zhang2018randomizedddd, qiliqun2021tttt}. In these works, 
some sketching tensors including the ones extracted from random sampling are formed. However, they have some limitations. For example, the Gaussian random tensor in \cite{zhang2018randomizedddd, qiliqun2021tttt} is defined as a tensor whose first frontal slice is created by the standard normal distribution and  other frontal slices are all zero; the random sampling tensor in \cite{tarzanagh2018fast, ma2021randomized,chen2021regularized,du2021randomized} is formed similarly, that is, its first frontal slice is a sampling matrix but other frontal slices are all zero. In this way, 
the transformed tensor by the discrete Fourier transform (DFT) along the third dimension 
will have the same frontal slices. On the other hand, 
a tensor problem based on t-product will be transformed into multiple independent matrix subproblems in the Fourier domain. Thus, 
the above sketching tensors 
will lead to the sketching matrices or the way for sampling in every matrix subproblem being the same. Taking the TRK method as an example, if we choose an index with the probabilities corresponding to the horizontal slices of $\mathcal{A}$, then every subsystem in the Fourier domain uses the same index to update at each iteration. Since these subsystems are independent, 
choosing different indices 
for different subsystems may be better. 

In \cite{ma2021randomized}, the authors also found the above limitation and mentioned that different indices can be selected for different subsystems. However, 
this strategy 
only works for complex-valued problems in the complex field. For real-valued problems in the field of real numbers, it is no longer feasible because 
the final solution is complex-valued. 
To the best of our knowledge, there is no work published to solve this problem in the real field. In this paper, we provide two improved strategies for our TSP method and its adaptive variants. The first one is based on an equivalence transformation, and the other is to take the real part of the last iterate directly. For the former, we present its theoretical guarantees. However,  
it is a little difficult to implement this method when combined with the adaptive sampling idea. For the latter, 
it has good performance in numerical experiments. However, we can't provide its theoretical guarantees at present. 

The paper is organized as follows. Section \ref{sec:pre} presents the notation and preliminaries. In Section \ref{sec:TSP}, we propose the TSP method and its adaptive variants. 
The implementation of the proposed methods in the Fourier domain is discussed in Section \ref{sec:TSPF}. In Section \ref{sec:TSPsub3}, we devise two improved strategies 
for the TSP method and its adaptive variants. The numerical results on synthetic and real data are provided in Section \ref{sec:exp}. Finally, we give the conclusion of the whole paper.


\section{Notation and preliminaries}\label{sec:pre}

Throughout this paper, scalars are denoted by lowercase letters, e.g., $x$; vectors are denoted by boldface lowercase letters, e.g., $\mathbf{x}$; matrices are denoted by boldface capital letters, e.g., $X$; higher-order tensors are denoted by Euler script letters, e.g., $\mathcal{X}$.

For a third-order tensor $\mathcal{X}$, its $(i,j,k)$-th element is represented by $\mathcal{X}_{(i,j,k)}$; its fiber 
is a one-dimensional array denoted by fixing two indices, e.g., $\mathcal{X}_{(:,j,k)}$, $\mathcal{X}_{(i,:,k)}$ and $\mathcal{X}_{(i,j,:)}$ respectively represent the $(j,k)$-th column, $(i,k)$-th row and $(i,j)$-th tube fiber; 
its slice 
is a two-dimensional array defined by fixing one index, e.g., $\mathcal{X}_{(i,:,:)}$, $\mathcal{X}_{(:,j,:)}$ and $\mathcal{X}_{(:,:,k)}$ respectively represent the $i$-th horizontal, $j$-th lateral and $k$-th frontal slice. 
For convenience, the frontal slice $\mathcal{X}_{(:,:,k)}$ is written as $\mathcal{X}_{(k)}$.

Before presenting the definition of t-product, we do some preparations. 
\begin{definition}[see \cite{kilmer2013third}]\rm
An element $\mathbf{x}\in \mathbb{R}^{1\times 1\times l}$ is called a tubal scalar of length $l$ and the set of all tubal scalars of length $l$ is denoted by $\mathbb{K}_{l}$; an element $\overrightarrow{\mathcal{X}} \in \mathbb{R}^{m\times 1\times l}$ 
is  called a vector of tubal scalars of length $l$ with size $m$ and the corresponding set is denoted by $\mathbb{K}^{m}_{l}$; 
an element $\mathcal{X} \in \mathbb{R}^{m\times n\times l}$ 
is  called a matrix of tubal scalars of length $l$ with size $m\times n$ and the corresponding set is denoted by $\mathbb{K}^{m\times n}_{l}$.
\end{definition}

Throughout this paper, we  will refer to tubal matrix and third-order tensor interchangeably.
For a tubal matrix $\mathcal{X}\in \mathbb{K}^{m\times n}_{l}$, as done in \cite{kilmer2013third, kilmer2011factorization}, define
$$\text{bcirc}(\mathcal{X})=
             \begin{bmatrix}
               \mathcal{X}_{(1)} & \mathcal{X}_{(l)} & \cdots & \mathcal{X}_{(2)}\\
               \mathcal{X}_{(2)} & \mathcal{X}_{(1)} & \cdots & \mathcal{X}_{(3)}  \\
              \vdots & \vdots & \ddots & \vdots \\
               \mathcal{X}_{(l)} & \mathcal{X}_{(l-1)} & \cdots & \mathcal{X}_{(1)} \\
            \end{bmatrix}
            ,~~\text{unfold}(\mathcal{X})=\begin{bmatrix}
                                           \mathcal{X}_{(1)} \\
                                           \mathcal{X}_{(2)} \\
                                           \vdots \\
                                           \mathcal{X}_{(l)}
                                         \end{bmatrix}
                                         ,\ \text{fold}(\text{unfold}(\mathcal{X}))=\mathcal{X}$$
and $\text{bcirc}^{-1}(\text{bcirc}(\mathcal{X}))=\mathcal{X}.$
\begin{definition}[t-product \cite{kilmer2011factorization}] \rm
Let $\mathcal{X}\in \mathbb{K}^{m\times n}_{l}$ and $\mathcal{Y}\in \mathbb{K}^{n\times p}_{l}$. Then the t-product $\mathcal{X}*\mathcal{Y}\in \mathbb{K}^{m\times p}_{l}$ is defined by
$$\mathcal{X}*\mathcal{Y}=\text{fold}(\text{bcirc}(\mathcal{X})\text{unfold}(\mathcal{Y})).$$
\end{definition}

Note that the matrix $\text{bcirc}(\mathcal{X})$ can be block diagonalized by the DFT matrix combined with the Kronecker product. Specifically, for a tubal matrix $\mathcal{X}\in \mathbb{K}^{m\times n}_{l}$ and the unitary DFT matrix $F_{l}\in \mathbb{C}^{l\times l}$,
\begin{align}\label{sec2e1}
\text{bdiag}(\widehat{\mathcal{X}})=(F_{l}\otimes I_{m})\text{bcirc}(\mathcal{X})(F_{l}^{H}\otimes I_{n})=\begin{bmatrix}
                                                                                 \widehat{\mathcal{ X}}_{(1)} &   &   &   \\
                                                                                    & \widehat{\mathcal{ X}}_{(2)} &   &   \\
                                                                                    &   & \ddots &   \\
                                                                                    &   &   &  \widehat{\mathcal{ X}}_{(l)} \\
                                                                               \end{bmatrix},
\end{align}
where $F_{l}^{H}$ denotes the conjugate transpose of $F_{l}$ and 
the matrices $\widehat{\mathcal{ X}}_{(k)}$ for $k=1,2,\cdots,l$ are the frontal slices of the tubal matrix $\widehat{\mathcal{ X}}$ which is obtained by applying the DFT on $\mathcal{ X}$ along the third dimension. 
We can use the Matlab function $\widehat{\mathcal{ X}}=\text{fft}(\mathcal{ X},[~],3)$ to  calculate $\widehat{\mathcal{ X}}$ directly, and use the inverse FFT to calculate $\mathcal{ X}$ from $\widehat{\mathcal{ X}}$, that is, $\mathcal{ X}=\text{ifft}(\widehat{\mathcal{ X}},[~],3)$. Thus, as noted in \cite{kilmer2011factorization}, the t-product $\mathcal{X}*\mathcal{Y}$ can be computed by computing FFT along each tubal fiber of $\mathcal{X}$ and $\mathcal{Y}$ to obtain $\widehat{\mathcal{ X}}=\text{fft}(\mathcal{ X},[~],3)$ and $\widehat{\mathcal{Y}}=\text{fft}(\mathcal{ Y},[~],3)$, multiplying each pair of the frontal slices of $\widehat{\mathcal{X}}$ and $\widehat{\mathcal{Y}}$ to get the frontal slices of $\widehat{\mathcal{Z}}$, and then taking inverse FFT along the third dimension of $\widehat{\mathcal{Z}}$ to get the desired result.

We also need the following definitions and the related results. 

\begin{definition}[transpose \cite{kilmer2011factorization}]\rm
For $\mathcal{X}\in \mathbb{K}^{m\times n}_{l}$, the transpose $\mathcal{X}^{T}$ is defined by taking the transpose of all the frontal slices and reversing the order of the second to last frontal slices.
\end{definition}

\begin{definition}[T-symmetric \cite{kilmer2011factorization}]\rm
For $\mathcal{X}\in \mathbb{K}^{n\times n}_{l}$, it is called T-symmetric if $\mathcal{X}=\mathcal{X}^{T}$. 
\end{definition}

\begin{definition}[identity tubal matrix \cite{kilmer2011factorization}]\rm
The identity tubal matrix $\mathcal{I}\in \mathbb{K}^{n\times n }_{l}$ is the tubal matrix whose first frontal slice is the $n\times n$ identity matrix and other frontal slices are all zero. 
\end{definition}

\begin{definition}[inverse \cite{kilmer2011factorization}]\rm
Let $\mathcal{X}\in \mathbb{K}^{n\times n}_{l}$. If there exists $\mathcal{Y} \in \mathbb{K}^{n\times n}_{l}$ such that
$$\mathcal{X}*\mathcal{Y}=\mathcal{I}~~~and~~~\mathcal{Y}*\mathcal{X}=\mathcal{I},$$
then $\mathcal{X}$ is said to be invertible, and $\mathcal{Y}$ is the inverse of $\mathcal{X}$, which is denoted by $\mathcal{X}^{-1}$.
\end{definition}

\begin{definition}[Moore-Penrose inverse \cite{jin2017generalized}]\rm
Let $\mathcal{X}\in \mathbb{K}^{m\times n}_{l}$. If there exists $\mathcal{Y} \in \mathbb{K}^{n\times m}_{l}$ such that
$$\mathcal{X}*\mathcal{Y}*\mathcal{X}=\mathcal{X},~~~\mathcal{Y}*\mathcal{X}*\mathcal{Y}=\mathcal{Y},~~~(\mathcal{X}*\mathcal{Y})^{T}=\mathcal{X}*\mathcal{Y},~~~(\mathcal{Y}*\mathcal{X})^{T}=\mathcal{Y}*\mathcal{X},$$
then $\mathcal{Y}$ is called the Moore-Penrose inverse of $\mathcal{X}$ and is denoted by $\mathcal{X}^{\dag}$.
\end{definition}


\begin{lemma}[see \cite{jin2017generalized}]
The Moore-Penrose inverse of any tubal matrix $\mathcal{X}\in \mathbb{K}^{m\times n}_{l}$ exists and is unique, and if $\mathcal{X}$ is invertible, then $\mathcal{X}^{\dag}=\mathcal{X}^{-1}$.
\end{lemma}

\begin{definition} [orthogonal tubal matrix \cite{kilmer2011factorization}]\rm
For $\mathcal{X}\in \mathbb{K}^{n\times n}_{l}$, it is orthogonal if $\mathcal{X}^{T}*\mathcal{X}=\mathcal{X}*\mathcal{X}^{T}=\mathcal{I}$.
\end{definition}

\begin{definition}[see \cite{kilmer2011factorization}]\rm
For 
$\mathcal{X} \in \mathbb{K}^{m\times n}_{l}$,  define
\begin{eqnarray*}
&\mathbf{Range}(\mathcal{X})=\left\{\overrightarrow{\mathcal{V}} \in \mathbb{K}^{m}_{l}|\overrightarrow{\mathcal{V}}=\mathcal{X}*\overrightarrow{\mathcal{Y}},~\text{for}~\text{any}~\overrightarrow{\mathcal{Y}}\in \mathbb{K}^{n}_{l}\right\},\
\mathbf{ Null}(\mathcal{X})=\left\{\overrightarrow{\mathcal{V}} \in \mathbb{K}^{n}_{l}|\mathcal{X}*\overrightarrow{\mathcal{V}}=\overrightarrow{\emph{O}}\right\},\\ 
&\mathbf{Colsp}(\mathcal{X})=\left\{\mathcal{V} \in \mathbb{K}^{m\times p}_{l}|~\text{for}~\text{all}~1\leq j\leq p,~\mathcal{V}_{(:,j,:)}\in \rm \mathbf{Range}(\mathcal{X})\right\}.
\end{eqnarray*}
\end{definition}

\begin{definition}[see \cite{kilmer2013third}]\rm
For 
$\mathcal{P}\in \mathbb{K}^{n\times n}_{l}$, it is a projector if $\mathcal{P}^{2}=\mathcal{P}*\mathcal{P}=\mathcal{P}$, and is orthogonal projector if $\mathcal{P}^{T}=\mathcal{P}$ also holds. 
\end{definition}

Note that $\mathcal{X}*(\mathcal{X}^{T}*\mathcal{X})^{\dag}*\mathcal{X}^{T}$ is an orthogonal projector onto $\mathbf{Range}(\mathcal{X})$.

\begin{definition}[T-symmetric T-positive (semi)definite \cite{zheng2021t}]\rm
For $\mathcal{X}\in \mathbb{K}^{n\times n }_{l}$, 
it is called T-symmetric T-positive (semi)definite if and only if $\mathcal{X}$ is  T-symmetric and $\langle\overrightarrow{\mathcal{Y}},\mathcal{X}*\overrightarrow{\mathcal{Y}}\rangle>(\geq) 0$ holds for any nonzero $\overrightarrow{\mathcal{Y}}\in \mathbb{K}^{n}_{l}$ (for any $\overrightarrow{\mathcal{Y}}\in \mathbb{K}^{n}_{l}$).
\end{definition}

\begin{proposition}[see \cite{qi2021t, zheng2021t}]
For 
$\mathcal{X}\in \mathbb{K}^{n\times n}_{l}$, it is T-symmetric if and only if $\text{bcirc}(\mathcal{X})$ is symmetric, 
is invertible if and only if $\text{bcirc}(\mathcal{X})$ is invertible, 
is orthogonal if and only if $\text{bcirc}(\mathcal{X})$ is orthogonal, and is T-symmetric T-positive (semi)definite if and only if $\text{bcirc}(\mathcal{X})$ is symmetric positive (semi)definite 
if and only if $\widehat{\mathcal{X}}_{(k)}$ for $k=1,2,\cdots,l$ are all Hermitian positive (semi)definite.
\end{proposition}

\begin{lemma}
Assume that $\mathcal{X}\in \mathbb{K}^{n\times n }_{l}$ is a T-symmetric T-positive (semi)definite tubal matrix, and define $\mathcal{X}^{\frac{1}{2}}=bcirc^{-1}(bcirc(\mathcal{X})^{\frac{1}{2}})$. Then $\mathcal{X}=\mathcal{X}^{\frac{1}{2}}*\mathcal{X}^{\frac{1}{2}}$ and $bcirc(\mathcal{X}^{\frac{1}{2}})=bcirc(\mathcal{X})^{\frac{1}{2}}$.
\end{lemma}
\emph{Proof:} Note that  $\text{bcirc}(\mathcal{A}*\mathcal{B})=\text{bcirc}(\mathcal{A})\text{bcirc}(\mathcal{B})$ holds for any $\mathcal{A}\in \mathbb{K}^{m\times n}_{l}$ and $\mathcal{B}\in \mathbb{K}^{n\times p}_{l}$, which can be found in  \cite{lund2020tensor}. Then, we can obtain $$\text{bcirc}^{-1}(\text{bcirc}(\mathcal{A})\text{bcirc}(\mathcal{B}))=\text{bcirc}^{-1}(\text{bcirc}(\mathcal{A}*\mathcal{B}))=\mathcal{A}*\mathcal{B}=\text{bcirc}^{-1}(\text{bcirc}(\mathcal{A}))*\text{bcirc}^{-1}(\text{bcirc}(\mathcal{B})).$$
Thus, considering $\mathcal{X}^{\frac{1}{2}}=\text{bcirc}^{-1}(\text{bcirc}(\mathcal{X})^{\frac{1}{2}})$, we have
\begin{align*}
  \mathcal{X}&=\text{bcirc}^{-1}(\text{bcirc}(\mathcal{X}))
  =\text{bcirc}^{-1}(\text{bcirc}(\mathcal{X})^{\frac{1}{2}}\text{bcirc}(\mathcal{X})^{\frac{1}{2}})\\
   &=\text{bcirc}^{-1}(\text{bcirc}(\mathcal{X})^{\frac{1}{2}})*\text{bcirc}^{-1}(\text{bcirc}(\mathcal{X})^{\frac{1}{2}})
    =\mathcal{X}^{\frac{1}{2}}*\mathcal{X}^{\frac{1}{2}},
 \end{align*}
 and
  \begin{align*}
  \text{bcirc}(\mathcal{X}^{\frac{1}{2}})&=\text{bcirc}(\text{bcirc}^{-1}(\text{bcirc}(\mathcal{X})^{\frac{1}{2}}))
  =\text{bcirc}(\mathcal{X})^{\frac{1}{2}}.
 \end{align*}
Then, the desired results hold.

\begin{definition}\rm
Let $Q \in \mathbb{R}^{n\times n}$ be a symmetric 
positive definite matrix. For any vectors $\mathbf{x}$, $\mathbf{y}\in \mathbb{R}^{n}$, their weighted inner product and the weighted induced norm are defined as $$\langle \mathbf{x},\mathbf{y}\rangle_{Q}\overset{\text{def}}{=}\langle Q\mathbf{x},\mathbf{y}\rangle \textrm{ and } \|\mathbf{x}\|_{Q}\overset{\text{def}}{=}\sqrt{\langle \mathbf{x},\mathbf{x}\rangle_{Q}},$$ respectively. For any matrix $M \in \mathbb{R}^{n\times p}$, its  weighted 2-norm and weighted Frobenius norm are defined as $$\|M\|_{2(Q)}\overset{\text{def}}{=}\max\limits_{{\mathbf{x}\in \mathbb{R}^{p},\|\mathbf{x}\|_{Q}=1}}{\|M\mathbf{x}\|_{Q}} \textrm{ and } \|M\|_{F(Q)}\overset{\text{def}}{=}\sqrt{\sum\limits_{j=1}^{p}\|M_{(:,j)}\|^{2}_{Q}},$$ respectively.
\end{definition}

Next, we extend the weighted norms for vectors and matrices to tubal vectors and tubal matrices, respectively.
\begin{definition}\rm
Let $\mathcal{Q}\in \mathbb{K}^{n\times n}_{l}$ be a T-symmetric 
T-positive definite tubal matrix. For any tubal vectors $\overrightarrow{\mathcal{X}}$, $\overrightarrow{\mathcal{Y}}\in \mathbb{K}^{n}_{l}$, their weighted inner product and the weighted induced norm are defined as
$$\langle\overrightarrow{\mathcal{X}},\overrightarrow{\mathcal{Y}}\rangle _{\mathcal{Q}}\overset{\text{def}}{=}\langle\mathcal{Q}*\overrightarrow{\mathcal{X}},\overrightarrow{\mathcal{Y}}\rangle=\langle \text{bcirc}(\mathcal{Q})\text{unfold}(\overrightarrow{\mathcal{X}}),\text{unfold}(\overrightarrow{\mathcal{Y}})\rangle=\langle \text{unfold}(\overrightarrow{\mathcal{X}}),\text{unfold}(\overrightarrow{\mathcal{Y}})\rangle _{\text{bcirc}(\mathcal{Q})},$$
and $$\|\overrightarrow{\mathcal{X}}\|_{\mathcal{Q}}\overset{\text{def}}{=}\sqrt{\langle\overrightarrow{\mathcal{X}},\overrightarrow{\mathcal{X}}\rangle _{\mathcal{Q}}}=\sqrt{\langle \text{unfold}(\overrightarrow{\mathcal{X}}),\text{unfold}(\overrightarrow{\mathcal{X}})\rangle _{\text{bcirc}(\mathcal{Q})}}=\|\text{unfold}(\overrightarrow{\mathcal{X}})\|_{\text{bcirc}(\mathcal{Q})},
$$
respectively. For any tubal matrix $\mathcal{M} \in \mathbb{K}^{n\times p}_{l}$, its weighted 2-norm and weighted Frobenius norm are defined as $$\|\mathcal{M}\|_{2(\mathcal{Q})}\overset{\text{def}}{=}\max\limits_{\overrightarrow{\mathcal{X}}\in \mathbb{K}^{p}_{l},\|\overrightarrow{\mathcal{X}}\|_{\mathcal{Q}}=1}{\|\mathcal{M}*\overrightarrow{\mathcal{X}}\|_{\mathcal{Q}}} \textrm{ and } \|\mathcal{M}\|_{F(\mathcal{Q})}\overset{\text{def}}{=}\sqrt{\sum\limits_{j=1}^{p}\|\mathcal{M}_{(:,j,:)}\|^{2}_{\mathcal{Q}}},$$ respectively.
\end{definition}

It is clear that
 \begin{align*}
  \|\mathcal{M}\|_{2(\mathcal{Q})}&=\max\limits_{\overrightarrow{\mathcal{X}}\in \mathbb{K}^{p}_{l},\|\overrightarrow{\mathcal{X}}\|_{\mathcal{Q}}=1}{\|\mathcal{M}*\overrightarrow{\mathcal{X}}\|_{\mathcal{Q}}}
  =\mathop{\max}_{\overrightarrow{\mathcal{X}}\in \mathbb{K}^{p}_{l},\|\text{unfold}(\overrightarrow{\mathcal{X}})\|_{\text{bcirc}(\mathcal{Q})}=1}{\|\text{unfold}(\mathcal{M}*\overrightarrow{\mathcal{X}})\|_{\text{bcirc}(\mathcal{Q})}}\\
   &=\mathop{\max}_{\text{unfold}(\overrightarrow{\mathcal{X}})\in \mathbb{R}^{pl},\|\text{unfold}(\overrightarrow{\mathcal{X}})\|_{\text{bcirc}(\mathcal{Q})}=1}{\|\text{bcirc}(\mathcal{M})\text{unfold}(\overrightarrow{\mathcal{X}})\|_{\text{bcirc}(\mathcal{Q})}}
   =\|\text{bcirc}(\mathcal{M})\|_{2(\text{bcirc}(\mathcal{Q}))},\\
  \|\mathcal{M}\|_{F(\mathcal{Q})}&=\sqrt{\sum_{j=1}^{p}\|\mathcal{M}_{(:,j,:)}\|^{2}_{\mathcal{Q}}}
  =\sqrt{\sum_{j=1}^{p}\|\text{unfold}(\mathcal{M})_{(:,j)}\|^{2}_{\text{bcirc}(\mathcal{Q})}}
   =\|\text{unfold}(\mathcal{M})\|_{F(\text{bcirc}(\mathcal{Q}))}.
 \end{align*}

\begin{lemma}
 Let $\mathcal{Q}\in \mathbb{K}^{n\times n}_{l}$ be a T-symmetric 
T-positive definite tubal matrix. 
Then for any tubal matrix $\mathcal{M} \in \mathbb{K}^{n\times p}_{l}$,
$$ \|\mathcal{M}\|_{F(\mathcal{Q})}=\|\mathcal{Q}^{\frac{1}{2}}*\mathcal{M}\|_{F}.$$
\end{lemma}
\emph{Proof:}
The result can be concluded by the properties of t-product and the definition of the weighted norm. Specifically,
 \begin{align*}
  \|\mathcal{M}\|_{F(\mathcal{Q})}&=\sqrt{\sum_{j=1}^{p}\|\text{unfold}(\mathcal{M})_{(:,j)}\|^{2}_{\text{bcirc}(\mathcal{Q})}}\\
  &=\sqrt{\sum_{j=1}^{p}(\text{unfold}(\mathcal{M})_{(:,j)})^{T}\text{bcirc}(\mathcal{Q})\text{unfold}(\mathcal{M})_{(:,j)}}\\
  &=\sqrt{\sum_{j=1}^{p}(\text{unfold}(\mathcal{M})_{(:,j)})^{T}\text{bcirc}(\mathcal{Q})^{\frac{1}{2}}\text{bcirc}(\mathcal{Q})^{\frac{1}{2}}\text{unfold}(\mathcal{M})_{(:,j)}}\\
  &=\sqrt{\sum_{j=1}^{p}\|\text{bcirc}(Q)^{\frac{1}{2}}\text{unfold}(\mathcal{M})_{(:,j)}\|^{2}_{2}}=\|\text{bcirc}(\mathcal{Q})^{\frac{1}{2}}\text{unfold}(\mathcal{M})\|_{F}\\
   &=\|\text{bcirc}(\mathcal{Q}^{\frac{1}{2}})\text{unfold}(\mathcal{M})\|_{F}=\|\text{unfold}(\mathcal{Q}^{\frac{1}{2}}*\mathcal{M})\|_{F}=\|\mathcal{Q}^{\frac{1}{2}}*\mathcal{M}\|_{F}.
 \end{align*}

For a symmetric 
positive semidefinite matrix $D\in \mathbb{R}^{n\times n}$, we write the seminorm induced by $D$ as $\|\mathbf{x}\|_{D}=\sqrt{\langle \mathbf{x},D\mathbf{x}\rangle}$.
Similarly, we can also define the seminorm induced by a T-symmetric
T-positive semidefinite tubal matrix $\mathcal{D}\in \mathbb{K}^{n\times n}_{l}$ as
$$\|\overrightarrow{\mathcal{X}}\|_{\mathcal{D}}\overset{\text{def}}{=}\sqrt{\langle\overrightarrow{\mathcal{X}},\mathcal{D}*\overrightarrow{\mathcal{X}}\rangle}=\|\text{unfold}(\overrightarrow{\mathcal{X}})\|_{\text{bcirc}(\mathcal{D})},$$ 
and, for any tubal matrix $\mathcal{M} \in \mathbb{K}^{n\times p}_{l}$, define  
$$\|\mathcal{M}\|_{2(\mathcal{D})}\overset{\text{def}}{=}\max\limits_{\overrightarrow{\mathcal{X}}\in \mathbb{K}^{p}_{l},\|\overrightarrow{\mathcal{X}}\|_{\mathcal{D}}=1}{\|\mathcal{M}*\overrightarrow{\mathcal{X}}\|_{\mathcal{D}}} \textrm{ and } \|\mathcal{M}\|_{F(\mathcal{D})}\overset{\text{def}}{=}\sqrt{\sum\limits_{j=1}^{p}\|\mathcal{M}_{(:,j,:)}\|^{2}_{\mathcal{D}}}.$$ 
Moreover, we have $ \|\mathcal{M}\|_{F(\mathcal{D})}=\|\mathcal{D}^{\frac{1}{2}}*\mathcal{M}\|_{F}.$

In addition, define $\triangle_{q}\overset{\text{def}}{=}\{\mathbf{p}=(p_{1},p_{2},\cdots,p_{q}) \in \mathbb{R}^{q} | \sum_{i=1}^{q} p_{i}=1,~p_{i}\geq0 \textrm{ for } i=1,\cdots,q\}$.
If $x_{i}$ depends on an index $i=1,2,\cdots,q$, we denote $\mathbf{E}_{i\sim \mathbf{p}}[x_{i}]\overset{\text{def}}{=}\sum_{i=1}^{q}p_{i}x_{i}$, where $\mathbf{p}\in \triangle_{q}$ and $i\sim \mathbf{p}$ means that the index $i$ is sampled with the probability $p_{i}$.  

\section{TSP method and its adaptive variants}\label{sec:TSP}

We first present the TSP method, and then introduce the adaptive sampling idea into the TSP method. 

\subsection{TSP method}\label{sec:TSPsub1}

Similar to the MSP method, the TSP method is designed to pursue the next iterate ${\mathcal{X}^{t+1}} \in {{\mathbb{K}}^{n\times p}_{l}}$  which is the nearest point to $\mathcal{X}^{t}$ and at the same time satisfies a sketched version of the problem (\ref{sec1e1}), that is
\begin{equation}\label{sec3.1e1}
\mathcal{X}^{t+1}=\mathop{\arg\min}_{\mathcal{X} \in {\mathbb{K}^{n\times p}_{l}}}\|\mathcal{X}-\mathcal{X}^{t}\|^2_{F(\mathcal{Q})}~\text{subject}~\text{to}~ \mathcal{S}^{T}*\mathcal{A}*\mathcal{X}=\mathcal{S}^{T}* \mathcal{B},
\end{equation}
where $\mathcal{S}\in {\mathbb{K}^{m\times \tau}_{l}}$ is a sketching tubal matrix drawn in an independent and identical distributed (i.i.d.) fashion from a fixed distribution $\mathfrak{D}$ at each iteration, and $\mathcal{Q} \in \mathbb{K}^{n\times n}_{l}$ is a T-symmetric T-positive definite tubal matrix. The distribution $\mathfrak{D}$ and tubal matrix $\mathcal{Q}$ are the parameters of the method. 
Making use of the algebraic properties of t-product, we can get the explicit solution to (\ref{sec3.1e1}) as
\begin{equation}\label{sec3.1e2}
\mathcal{X}^{t+1}=\mathcal{X}^{t}-\mathcal{Q}^{-1}*\mathcal{A}^{T}*\mathcal{S}*(\mathcal{S}^{T}*\mathcal{A}*\mathcal{Q}^{-1}*\mathcal{A}^{T}*\mathcal{S})^{\dag}*\mathcal{S}^{T}*(\mathcal{A}*\mathcal{X}^{t}-\mathcal{B}),
\end{equation}
and then we obtain the TSP method, i.e., Algorithm \ref{TSP}.

\begin{algorithm}[htbp]
\caption{TSP method}\label{TSP}
$\textbf{Input:}$ $\mathcal{X}^{0}\in \mathbb{K}^{n\times p}_{l}$, $\mathcal{A}\in \mathbb{K}^{m\times n }_{l}$, $\mathcal{B}\in \mathbb{K}^{m\times p }_{l}$

$\textbf{Parameters:}$ fixed distribution $\mathfrak{D}$ over random tubal matrices, T-symmetric T-positive definite tubal matrix $\mathcal{Q}\in \mathbb{K}^{n\times n}_{l}$

$\textbf{for}$ $t=0,1,2,\cdots$

~~~~~Sample an independent copy $\mathcal{S}\sim \mathfrak{D}$

~~~~~Compute $\mathcal{G}=\mathcal{S}*(\mathcal{S}^{T}*\mathcal{A}*\mathcal{Q}^{-1}*\mathcal{A}^{T}*\mathcal{S})^{\dag}*\mathcal{S}^{T}$

~~~~~$\mathcal{X}^{t+1}=\mathcal{X}^{t}-\mathcal{Q}^{-1}*\mathcal{A}^{T}*\mathcal{G}*(\mathcal{A}*\mathcal{X}^{t}-\mathcal{B})$

$\textbf{end for}$

$\textbf{Output:}$ last iterate $\mathcal{X}^{t+1}$
\end{algorithm}

\begin{remark}\rm\label{sec3.1rem2.1}
If we choose 
$\mathcal{S}=\mathcal{I}_{(:,i,:)}\in \mathbb{K}^{m}_{l}$ with $i=1,2,\cdots,m$ 
and $\mathcal{Q}=\mathcal{I} \in \mathbb{K}^{n\times n }_{l}$, 
then it follows from (\ref{sec3.1e2}) that
\begin{align*} 
\mathcal{X}^{t+1}=\mathcal{X}^{t}-(\mathcal{A}_{(i,:,:)})^{T}*\left(\mathcal{A}_{(i,:,:)}*(\mathcal{A}_{(i,:,:)})^{T}\right)^{\dag}*\left(\mathcal{A}_{(i,:,:)}*\mathcal{X}^{t}-\mathcal{B}_{(i,:,:)}\right).
\end{align*}
When $i$ is selected uniformly or with the probabilities proportional to the norm of horizontal slices\footnote{In this case, the distribution $\mathfrak{D}$ is a discrete distribution and $\mathbf{P}(\mathcal{S}=\mathcal{I}_{(:,i,:)})=p_{i}$ with $i=1,\cdots,m$, where $p_{i}=\frac{1}{m}$ or $p_{i}=\frac{\|\mathcal{A}_{(i,:,:)}\|_{F}^{2}}{\|\mathcal{A}\|_{F}^{2}}$.}, the TSP method 
will reduce to the TRK method in \cite{ma2021randomized}.
\end{remark}

Next, we shall discuss the convergence analysis of the TSP method.

\begin{theorem}\label{sec3.1thm3.1}
With the notation in Algorithm \ref{TSP}, assume that 
$\mathbf{E}[\mathcal{Z}]$ is T-symmetric T-positive definite with probability $1$, where $\mathcal{Z}=\mathcal{Q}^{-\frac{1}{2}}*\mathcal{W}*Q^{-\frac{1}{2}}$ and $\mathcal{W}=\mathcal{A}^{T}*\mathcal{G}*\mathcal{A}$, $\mathcal{X}^{\star}$ satisfies $\mathcal{A}*\mathcal{X}^{\star}=\mathcal{B}$, and $\mathcal{X}^{t}$ is the $t$-th approximation of $\mathcal{X}^{\star}$ 
with initial iterate $\mathcal{X}^{0}$. Then
\begin{align}
\mathbf{E}\left[\|\mathcal{X}^{t}-\mathcal{X}^{\star}\|_{F(\mathcal{Q})}^{2}|\mathcal{X}^{0}\right]\leq\left(1-\lambda_{\mathop{\min}}(\mathbf{E}[\rm{bcirc}(\mathcal{Z})])\right)^{t}\|\mathcal{X}^{0}-\mathcal{X}^{\star}\|_{F(\mathcal{Q})}^{2}.
\end{align}
\end{theorem}
\emph{Proof:} Combining (\ref{sec3.1e2}) and the fact that $\mathcal{A}*\mathcal{X}^{\star}=\mathcal{B}$, we have
\begin{align}
\mathcal{X}^{t+1}-\mathcal{X}^{\star}=(\mathcal{I}-\mathcal{Q}^{-1}*\mathcal{W})*(\mathcal{X}^{t}-\mathcal{X}^{\star}).\label{sec3.1e4}
\end{align}
Multiplying both sides of (\ref{sec3.1e4}) by $\mathcal{Q}^{\frac{1}{2}}$, we obtain
\begin{align*}
\mathcal{Q}^{\frac{1}{2}}*(\mathcal{X}^{t+1}-\mathcal{X}^{\star})&=\mathcal{Q}^{\frac{1}{2}}*(\mathcal{I}-\mathcal{Q}^{-1}*\mathcal{W})*\mathcal{Q}^{-\frac{1}{2}}*\mathcal{Q}^{\frac{1}{2}}*(\mathcal{X}^{t}-\mathcal{X}^{\star})
=(\mathcal{I}-\mathcal{Z})*\mathcal{Q}^{\frac{1}{2}}*(\mathcal{X}^{t}-\mathcal{X}^{\star}).
\end{align*}
Let $\Gamma^{t}=\mathcal{Q}^{\frac{1}{2}}*(\mathcal{X}^{t}-\mathcal{X}^{\star})$. Thus, the above equation can be rewritten as
$\Gamma^{t+1}=(\mathcal{I}-\mathcal{Z})*\Gamma^{t}.$
Applying the Frobenius norm to its two sides, we get
\begin{align}\label{sec3.1e6}
\|\Gamma^{t+1}\|_{F}^{2}=\|(\mathcal{I}-\mathcal{Z})*\Gamma^{t}\|^{2}_{F}=\|\Gamma^{t}\|^{2}_{F}-\|\mathcal{Z}*\Gamma^{t}\|^{2}_{F},
\end{align}
where the second equality is from 
the Pythagorean theorem. By taking expectation conditioned on $\mathcal{X}^{t}$, we have
\begin{align}\label{sec3.1e7}
\mathbf{E}[\|\Gamma^{t+1}\|_{F}^{2}|\mathcal{X}^{t}]
  &=\|\Gamma^{t}\|_{F}^{2}-\mathbf{E}[\|\mathcal{Z}*\Gamma^{t}\|^{2}_{F}].
 \end{align}
Note that
\begin{align*}
\mathbf{E}[\|\mathcal{Z}*\Gamma^{t}\|^{2}_{F}]&=\mathbf{E}[\|\text{bcirc}(\mathcal{Z})\text{unfold}(\Gamma^{t})\|^{2}_{F}]=\sum_{j=1}^{p}\langle\mathbf{E}[\text{bcirc}(\mathcal{Z})]\text{unfold}(\Gamma^{t})_{(:,j)}, \text{unfold}(\Gamma^{t})_{(:,j)}\rangle\\
&\geq\lambda_{\min}(\mathbf{E}[\text{bcirc}(\mathcal{Z})])\sum_{j=1}^{p}\|\text{unfold}(\Gamma^{t})_{(:,j)}\|^{2}_{2}=\lambda_{\min}(\mathbf{E}[\text{bcirc}(\mathcal{Z})])\|\Gamma^{t}\|^{2}_{F},
\end{align*}
where the inequality follows from the assumption that $\mathbf{E}[\mathcal{Z}]$ is T-symmetric T-positive definite with probability $1$. Therefore,
\begin{align*}
\mathbf{E}[\|\Gamma^{t+1}\|_{F}^{2}|\mathcal{X}^{t}]
  &\leq (1-\lambda_{\min}(\mathbf{E}[\text{bcirc}(\mathcal{Z})]))\|\Gamma^{t}\|_{F}^{2},
 \end{align*}
 that is
\begin{align*}
  \mathbf{E}[\|\mathcal{X}^{t+1}-\mathcal{X}^{\star}\|_{F(\mathcal{Q})}^{2}|\mathcal{X}^{t}] \leq (1-\lambda_{\min}(\mathbf{E}[ \text{bcirc}(\mathcal{Z})])\|\mathcal{X}^{t}-\mathcal{X}^{\star}\|_{F(\mathcal{Q})}^{2}.
 \end{align*}
 Taking expectation again and unrolling the recurrence give the desired result.

\begin{remark}\rm\label{sec3.1rem3.1}
Now, we show that the convergence rate $\rho_{TSP}=1-\lambda_{\min}(\mathbf{E}[ \text{bcirc}(\mathcal{Z})])$ is smaller than 1. It is easy to check that $\text{bcirc}(\mathcal{Z})$ is an orthogonal projection and hence has eigenvalues $0$ or $1$. Furthermore, it 
projects onto a $d$-dimensional subspace $$\mathbf{Range}(\text{bcirc}(\mathcal{Q})^{-\frac{1}{2}}\text{bcirc}(\mathcal{A})^{T}\text{bcirc}(\mathcal{S})),$$ 
where $d\overset{def}=\mathbf{Rank}(\text{bcirc}(\mathcal{S})^{T}\text{bcirc}(\mathcal{A})) \leq \min(ql,nl)$. Using Jensen's inequality, as well as the fact that both $A \longmapsto \lambda_{\max}(A)$ and $A \longmapsto -\lambda_{\min}(A)$ are convex on the symmetric matrices, we can conclude that the spectrum of $\mathbf{E}[\text{bcirc}(\mathcal{Z})]$ is contained in $[0,1]$. Next, we turn to refine the lower and upper bounds of $\rho_{TSP}$. It follows that
\begin{align*}
\mathbf{E}[d]&=\mathbf{E}[\mathbf{Tr}( \text{bcirc}(\mathcal{Z}))]
=\mathbf{Tr}(\mathbf{E}[\text{bcirc}(\mathcal{Z})]) \geq nl\lambda_{\min}(\mathbf{E}[\text{bcirc}(\mathcal{Z})]),
\end{align*}
where the inequality holds because the trace of a matrix is equal to the sum of its eigenvalues. Thus, we have $\lambda_{\min}(\mathbf{E}[\text{bcirc}(\mathcal{Z})])\leq\frac{E[d]}{nl}$. Furthermore, $\mathbf{E}[\text{bcirc}(\mathcal{Z})]$ is symmetric positive definite because $\mathbf{E}[\mathcal{Z}]$ is T-symmetric T-positive definite, which immediately yields $\lambda_{\min}(\mathbf{E}[\text{bcirc}(\mathcal{Z})])> 0$.
All together, we have the following lower and upper bounds on $\rho_{TSP}$:
$$0\leq 1-\frac{\mathbf{E}[d]}{nl}\leq\rho_{TSP}<1.$$
So, the rate is indeed smaller than 1 and hence the sequence $\{\mathcal{X}^{t}\}_{t=0}^{\infty}$ generated by the TSP method can converge to $\mathcal{X}^{\star}$.
\end{remark}

\begin{remark}\rm\label{sec3.1rem2}
The convergence guarantee for the TRK method presented in \cite{ma2021randomized} is a special case of Theorem \ref{sec3.1thm3.1}. Specifically, choosing 
$\mathcal{S}=\mathcal{I}_{(:,i,:)}$ with $i=1,2,\cdots,m$  and $\mathcal{Q}=\mathcal{I} \in \mathbb{K}^{n\times n }_{l}$ in Theorem \ref{sec3.1thm3.1}, we can recover the result given in Theorem 3.1 in \cite{ma2021randomized}.
\end{remark}

\subsection{Three adaptive TSP methods}\label{sec:ATSP}
As shown in Section \ref{sec:TSPsub1}, a key step of the TSP method is to choose a sketching tubal matrix $\mathcal{S}$ in an i.i.d. fashion from a fixed distribution. In this subsection, we mainly study the 
adaptive sampling strategies on a finite set of sketching tubal matrices which is selected from a certain distribution in advance (the selection of the finite set is not considered in this paper). That is, letting $\boldsymbol{\mathcal{S}} = \{ \mathcal{S}_{i} \in \mathbb{K}^{m \times \tau}_{l},~\text{for}~ i=1,\cdots,q,~q \in \mathbb{N}\}$ be a finite set of sketching tubal matrices where $\tau\in \mathbb{N}$ is the sketch size, we want to choose $\mathcal{S}=\mathcal{S}_{i}$ from $\boldsymbol{\mathcal{S}}$ using adaptive sampling strategies. 
If the sampling probability distribution at each iteration is fixed, we call the corresponding method 
the nonadaptive TSP (NTSP) method, which is summarized in Algorithm \ref{NTSP}.
\begin{algorithm}[htbp]
\caption{NTSP method}\label{NTSP}
$\textbf{Input:}$ $\mathcal{X}^{0}\in \mathbb{K}^{n\times p}_{l}$, $\mathcal{A}\in \mathbb{K}^{m\times n}_{l}$, $\mathcal{B}\in \mathbb{K}^{m\times p}_{l}$, and $\mathbf{p}\in \triangle_{q}$

$\textbf{Parameters:}$ a set of sketching tubal matrices $\boldsymbol{\mathcal{S}}=[\mathcal{S}_{1},\cdots,\mathcal{S}_{q}]$, T-symmetric T-positive definite tubal matrix $\mathcal{Q}\in \mathbb{K}^{n\times n}_{l}$

$\textbf{for}$ $t=0,1,2,\cdots$ 

~~~~~$i^{t}\sim \mathbf{p}$

~~~~~Compute $\mathcal{G}_{i^{t}}=\mathcal{S}_{i^{t}}*(\mathcal{S}_{i^{t}}^{T}*\mathcal{A}*\mathcal{Q}^{-1}*\mathcal{A}^{T}*\mathcal{S}_{i^{t}})^{\dag}*\mathcal{S}_{i^{t}}^{T}$

~~~~~$\mathcal{X}^{t+1}=\mathcal{X}^{t}-\mathcal{Q}^{-1}*\mathcal{A}^{T}*\mathcal{G}_{i^{t}}*(\mathcal{A}*\mathcal{X}^{t}-\mathcal{B})$

$\textbf{end for}$

$\textbf{Output:}$ last iterate $\mathcal{X}^{t+1}$
\end{algorithm}

\begin{remark}\rm
There are some subtle differences between Algorithms \ref{TSP} and \ref{NTSP}. Specifically, the former draws a sketching tubal matrix from a fixed distribution at each iteration, while the latter needs to select a finite set of sketching tubal matrices from a distribution in advance, and then picks one from the finite set  with a fixed probability at each iteration. It should be noted that when the TSP method reduces to the TRK method, the two algorithms are the same.
\end{remark}

\subsubsection{Three adaptive sampling strategies}\label {sec:ASTP-AS}
Considering that the fixed sampling strategy in Algorithm \ref{NTSP} may choose a terrible $\mathcal{S}_{i}$ and hence leads to a bad convergence, 
we introduce three adaptive sampling strategies which use information about the current iterate.


Specifically, setting $\mathcal{S}=\mathcal{S}_{i^{t}}$ in (\ref{sec3.1e6}) and using the fact that $\mathcal{Z}_{i^{t}}=\mathcal{Q}^{-\frac{1}{2}}*\mathcal{A}^{T}*\mathcal{G}_{i^{t}}*\mathcal{A}*\mathcal{Q}^{-\frac{1}{2}}$ is an orthogonal projector onto $\mathbf{Range}(\mathcal{Q}^{-\frac{1}{2}}*\mathcal{A}^{T}*\mathcal{S}_{i^{t}})$, we have
\begin{align}
\|\mathcal{X}^{t+1}-\mathcal{X}^{\star}\|_{F(\mathcal{Q})}^{2}
&=\|\mathcal{X}^{t}-\mathcal{X}^{\star}\|_{F(\mathcal{Q})}^{2}-\|\mathcal{Z}_{i^{t}}*\mathcal{Q}^{\frac{1}{2}}*(\mathcal{X}^{t}-\mathcal{X}^{\star})\|_{F}^{2}\nonumber\\
&=\|\mathcal{X}^{t}-\mathcal{X}^{\star}\|_{F(\mathcal{Q})}^{2}-\|\mathcal{Q}^{\frac{1}{2}}*(\mathcal{X}^{t}-\mathcal{X}^{\star})\|_{F(\mathcal{Z}_{i^{t}})}^{2}\nonumber\\
&=\|\mathcal{X}^{t}-\mathcal{X}^{\star}\|_{F(\mathcal{Q})}^{2}-f_{i^{t}}(\mathcal{X}^{t}),\label{sec4.1e1}
\end{align}
which shows that the quantity of the error $\mathcal{X}^{t+1}-\mathcal{X}^{\star}$ is determined by $f_{i^{t}}(\mathcal{X}^{t})$. Consequently, in order to make the most progress in one step, we should choose $i^{t}$ corresponding to the largest sketched loss $f_{i^{t}}(\mathcal{X}^{t})$.
Since $\mathcal{X}^{\star}$ is unknown in practice,  we first rewrite $f_{i^{t}}(\mathcal{X}^{t})$ as
\begin{align*}
f_{i^{t}}(\mathcal{X}^{t})
&=\|\mathcal{Q}^{\frac{1}{2}}*(\mathcal{X}^{t}-\mathcal{X}^{\star})\|_{F(\mathcal{Z}_{i^{t}})}^{2}
=\sum_{j=1}^{p}\|\text{unfold}(\mathcal{Q}^{\frac{1}{2}}*(\mathcal{X}^{t}-\mathcal{X}^{\star}))_{(:,j)}\|_{\text{bcirc}(\mathcal{Z}_{i^{t}})}^{2}\\
&=\sum_{j=1}^{p}(\text{unfold}(\mathcal{X}^{t}-\mathcal{X}^{\star})_{(:,j)})^{T}\text{bcirc}(\mathcal{Q})^{\frac{1}{2}}\text{bcirc}(\mathcal{Z}_{i^{t}})\text{bcirc}(\mathcal{Q})^{\frac{1}{2}}\text{unfold}(\mathcal{X}^{t}-\mathcal{X}^{\star})_{(:,j)}\\
&=\sum_{j=1}^{p}(\text{unfold}(\mathcal{X}^{t}-\mathcal{X}^{\star})_{(:,j)})^{T}\text{bcirc}(\mathcal{A})^{T}\text{bcirc}(\mathcal{G}_{i^{t}})\text{bcirc}(\mathcal{A})\text{unfold}(\mathcal{X}^{t}-\mathcal{X}^{\star})_{(:,j)}\\
&=\sum_{j=1}^{p}\|\text{bcirc}(\mathcal{A})\text{unfold}(\mathcal{X}^{t}-\mathcal{X}^{\star})_{(:,j)}\|_{\text{bcirc}(\mathcal{G}_{i^{t}})}^{2}
=\|\mathcal{A}*(\mathcal{X}^{t}-\mathcal{X}^{\star})\|_{F(\mathcal{G}_{i^{t}})}^{2}
=\|\mathcal{A}*\mathcal{X}^{t}-\mathcal{B}\|_{F(\mathcal{G}_{i^{t}})}^{2}.
\end{align*}
Thus, according to (\ref{sec4.1e1}), we can present the first adaptive sampling strategy as follows
\begin{align}\label{sec4.1e3}
i^{t}=\mathop{\arg\max}_{i=1,\cdots,q}f_{i}(\mathcal{X}^{t})=\mathop{\arg\max}_{i=1,\cdots,q}\|\mathcal{A}*\mathcal{X}^{t}-\mathcal{B}\|_{F(\mathcal{G}_{i})}^{2},
\end{align}
which can be called the max-distance selection rule and the corresponding algorithm is described in Algorithm \ref{ATSP-MD}. 

\begin{algorithm}[htbp]
\caption{ATSP-MD method}\label{ATSP-MD}
$\textbf{Input:}$ $\mathcal{X}^{0}\in \mathbb{K}^{n\times p}_{l}$, $\mathcal{A}\in \mathbb{K}^{m\times n }_{l}$, and $\mathcal{B}\in \mathbb{K}^{m\times p }_{l}$

$\textbf{Parameters:}$ a set of sketching tubal matrices $\boldsymbol{\mathcal{S}}=[\mathcal{S}_{1},\cdots,\mathcal{S}_{q}]$, T-symmetric T-positive definite tubal matrix $\mathcal{Q}\in \mathbb{K}^{n\times n}_{l}$

$\textbf{for}$ $t=0,1,2,\cdots$

~~~~~$f_{i}(\mathcal{X}^{t})=\|\mathcal{A}*\mathcal{X}^{t}-\mathcal{B}\|_{F(\mathcal{G}_{i})}^{2}$ for $i=1,\cdots,q$

~~~~~$i^{t}=\mathop{\arg\max}\limits_{i=1,\cdots,q}f_{i}(\mathcal{X}^{t})$

~~~~~Compute $\mathcal{G}_{i^{t}}=\mathcal{S}_{i^{t}}*\left(\mathcal{S}_{i^{t}}^{T}*\mathcal{A}*\mathcal{Q}^{-1}*\mathcal{A}^{T}*\mathcal{S}_{i^{t}}\right)^{\dag}*\mathcal{S}_{i^{t}}^{T}$

~~~~~$\mathcal{X}^{t+1}=\mathcal{X}^{t}-\mathcal{Q}^{-1}*\mathcal{A}^{T}*\mathcal{G}_{i^{t}}*(\mathcal{A}*\mathcal{X}^{t}-\mathcal{B})$

$\textbf{end for}$

$\textbf{Output:}$ last iterate $\mathcal{X}^{t+1}$
\end{algorithm}

Now, we consider the expected decrease of the error $\mathcal{X}^{t+1}-\mathcal{X}^{\star}$. Let $\mathbf{p}^{t}\in \triangle_{q}$ and $i^{t}\sim \mathbf{p}^{t}$, where  $\mathbf{p}^{t}\overset{\text{def}}{=}(p^{t}_{1},\cdots,p^{t}_{q})$ with  $p_{i}^{t}=\mathbf{P}[\mathcal{S}_{i^{t}}=\mathcal{S}_{i}|\mathcal{X}^{t}]$ for $i=1,\cdots,q$, i.e., $p_{i}^{t}$ is the probability of $\mathcal{S}_{i}$ being sampled at the $t$-th iteration. Taking expectation conditioned on $\mathcal{X}^{t}$ in (\ref{sec4.1e1}), we have
\begin{align}\label{sec4.1e2}
\mathbf{E}[\|\mathcal{X}^{t+1}-\mathcal{X}^{\star}\|_{F(\mathcal{Q})}^{2}|\mathcal{X}^{t}]=\|\mathcal{X}^{t}-\mathcal{X}^{\star}\|_{F(\mathcal{Q})}^{2}-\mathbf{E}_{i\sim \mathbf{p}^{t}}[f_{i}(\mathcal{X}^{t})],
\end{align}
which tells us that if we want $\mathbf{E}[\|\mathcal{X}^{t+1}-\mathcal{X}^{\star}\|_{F(\mathcal{Q})}^{2}|\mathcal{X}^{t}]$ to be as small as possible, we should choose adaptive probabilities to make $\mathbf{E}_{i\sim \mathbf{p}^{t}}[f_{i}(\mathcal{X}^{t})]$ as large as possible. Since $\mathbf{E}_{i\sim \mathbf{p}^{t}}[f_{i}(\mathcal{X}^{t})]=\sum_{i=1}^{q}p_{i}^{t}f_{i}(\mathcal{X}^{t})$, we can achieve the above goal by sampling the indices corresponding to larger sketched losses with higher probability. An intuitive way is to choose the probabilities proportional to the sketched losses and we refer to such strategy as adaptive probabilities rule. The algorithm is summarized in Algorithm \ref{ATSP-PR}. 

\begin{algorithm}[htbp]
\caption{ATSP-PR method}\label{ATSP-PR}
$\textbf{Input:}$ $\mathcal{X}^{0}\in \mathbb{K}^{n\times p  }_{l}$, $\mathcal{A}\in \mathbb{K}^{m\times n }_{l}$, and $\mathcal{B}\in \mathbb{K}^{m\times p }_{l}$

$\textbf{Parameters:}$ a set of sketching tubal matrices $\boldsymbol{\mathcal{S}}=[\mathcal{S}_{1},\cdots,\mathcal{S}_{q}]$, T-symmetric T-positive definite tubal matrix $\mathcal{Q}\in \mathbb{K}^{n\times n }_{l}$

$\textbf{for}$ $t=0,1,2,\cdots$

~~~~~$f_{i}(\mathcal{X}^{t})=\|\mathcal{A}*\mathcal{X}^{t}-\mathcal{B}\|_{F(\mathcal{G}_{i})}^{2}$ for $i=1,\cdots,q$

~~~~~Calculate $\mathbf{p}^{t}\in \triangle_{q}$ such that $p_{i}^{t}=f_{i}(\mathcal{X}^{t})/(\sum_{i=1}^{q}f_{i}(\mathcal{X}^{t}))$ for $i=1,\cdots,q$

~~~~~$i^{t}\sim \mathbf{p}^{t}$

~~~~~Compute $\mathcal{G}_{i^{t}}=\mathcal{S}_{i^{t}}*\left(\mathcal{S}_{i^{t}}^{T}*\mathcal{A}*\mathcal{Q}^{-1}*\mathcal{A}^{T}*\mathcal{S}_{i^{t}}\right)^{\dag}*\mathcal{S}_{i^{t}}^{T}$

~~~~~$\mathcal{X}^{t+1}=\mathcal{X}^{t}-\mathcal{Q}^{-1}*\mathcal{A}^{T}*\mathcal{G}_{i^{t}}*(\mathcal{A}*\mathcal{X}^{t}-\mathcal{B})$

$\textbf{end for}$

$\textbf{Output:}$ last iterate $\mathcal{X}^{t+1}$
\end{algorithm}

In addition, there is another effective strategy, which aims to capture the indices corresponding to larger sketched losses as far as possible at each iteration. To this end, it considers removing the indices corresponding to the smaller sketched losses. 
To be specific, we first define an index set 
\begin{align}\label{sec2.e11111}
\mathfrak{W}_{t}=\left\{i|f_{i}(\mathcal{X}^{t})\geq\theta\mathop{\max}_{j=1,2,\cdots,q}f_{j}(\mathcal{X}^{t})+(1-\theta)\mathbf{E}_{j\sim \mathbf{p}}[f_{j}(\mathcal{X}^{t})]\right\},
\end{align}
where $\mathbf{p}\in \triangle_{q}$ and $\theta\in [0,1]$. Then, we choose the probabilities $\mathbf{p}^{t}\in \triangle_{q}$ such that 
\begin{equation}\label{sec2.e22222}
p_{i}^{t}=\left\{
\begin{array}{lcl}
\frac{f_{i}(\mathcal{X}^{t})}{\sum_{i\in \mathfrak{W}_{t}}f_{i}(\mathcal{X}^{t})}& & i\in \mathfrak{W}_{t}\\
0 & & i\notin \mathfrak{W}_{t}.
\end{array} \right.
\end{equation}
We call this strategy the capped sampling rule, which is summarized in Algorithm \ref{ATSP-CS}.

\begin{algorithm}[htbp]
\caption{ATSP-CS method}\label{ATSP-CS}
$\textbf{Input:}$ $\mathcal{X}^{0}\in \mathbb{K}^{n\times p }_{l}$, $\mathcal{A}\in \mathbb{K}^{m\times n }_{l}$, $\mathcal{B}\in \mathbb{K}^{m\times p }_{l}$, $\mathbf{p}\in \triangle_{q}$, and $\theta\in [0,1]$

$\textbf{Parameters:}$ a set of sketching tubal matrices $\boldsymbol{\mathcal{S}}=[\mathcal{S}_{1},\cdots,\mathcal{S}_{q}]$, T-symmetric T-positive definite tubal matrix $\mathcal{Q}\in \mathbb{K}^{n\times n }_{l}$

$\textbf{for}$ $t=0,1,2,\cdots$

~~~~~$f_{i}(\mathcal{X}^{t})=\|\mathcal{A}*\mathcal{X}^{t}-\mathcal{B}\|_{F(\mathcal{G}_{i})}^{2}$ for $i=1,\cdots,q$

~~~~~Determine the index set $\mathfrak{W}_{t}$, which is defined in (\ref{sec2.e11111})

~~~~~Calculate $\mathbf{p}^{t}\in \triangle_{q}$, which is defined in (\ref{sec2.e22222})

~~~~~$i^{t}\sim \mathbf{p}^{t}$

~~~~~Compute $\mathcal{G}_{i^{t}}=\mathcal{S}_{i^{t}}*\left(\mathcal{S}_{i^{t}}^{T}*\mathcal{A}*\mathcal{Q}^{-1}*\mathcal{A}^{T}*\mathcal{S}_{i^{t}}\right)^{\dag}*\mathcal{S}_{i^{t}}^{T}$

~~~~~$\mathcal{X}^{t+1}=\mathcal{X}^{t}-\mathcal{Q}^{-1}*\mathcal{A}^{T}*\mathcal{G}_{i^{t}}*(\mathcal{A}*\mathcal{X}^{t}-\mathcal{B})$

$\textbf{end for}$

$\textbf{Output:}$ last iterate $\mathcal{X}^{t+1}$
\end{algorithm}

\begin{remark}\rm
In the TRK setting, the above ATSP-MD, ATSP-PR and ATSP-CS methods are typically referred to as the ATRK-MD, ATRK-PR and ATRK-CS methods, which are the tensor versions of the greedy or adaptive MRK methods given in \cite{nutini2016convergence, Bai2018, Bai2018r}.
\end{remark}

\subsubsection{Convergence}\label{sec:ASTP-CON}
In this subsection, we discuss the convergence analysis of the nonadaptive and adaptive TSP methods proposed above. Before the formal discussions, we first prove two lemmas.
\begin{lemma}
With the notation in the NTSP, ATSP-MD, ATSP-PR, and ATSP-CS methods, let $\mathbf{p}\in \triangle_{q}$ and define
\begin{align}
\delta_{\infty}^{2}(\mathcal{Q},\boldsymbol{\mathcal{S}})\overset{def}{=}\mathop{\min}_{\overrightarrow{\mathcal{V}}\in \mathbf{Range}(\mathcal{Q}^{-1}*\mathcal{A}^{T})}\mathop{\max}_{i=1,\cdots,q}\frac{\|\mathcal{Q}^{\frac{1}{2}}*\overrightarrow{\mathcal{V}}\|_{\mathcal{Z}_{i}}^{2}}{\|\overrightarrow{\mathcal{V}}\|_{\mathcal{Q}}^{2}},\label{sec4.2e1}\\
\delta_{\mathbf{p}}^{2}(\mathcal{Q},\boldsymbol{\mathcal{S}})\overset{def}{=}\mathop{\min}_{\overrightarrow{\mathcal{V}}\in \mathbf{Range}(\mathcal{Q}^{-1}*\mathcal{A}^{T})}\frac{\|\mathcal{Q}^{\frac{1}{2}}*\overrightarrow{\mathcal{V}}\|_{\mathbf{E}_{i\sim \mathbf{p}}[\mathcal{Z}_{i}]}^{2}}{\|\overrightarrow{\mathcal{V}}\|_{\mathcal{Q}}^{2}},\label{sec4.2e2}
\end{align}
where $\mathcal{Z}_{i}$ is the same as $\mathcal{Z}_{i^t}$ defined above except that $i^t$ is replaced by $i$. 
Let $\mathcal{X}^{\star}$ satisfy $\mathcal{A}*\mathcal{X}^{\star}=\mathcal{B}$ and $\mathcal{X}^{t}$ be the $t$-th approximation of $\mathcal{X}^{\star}$ calculated by any nonadaptive and adaptive algorithms with initial iterate $\mathcal{X}^{0}\in \mathbf{Colsp}(\mathcal{Q}^{-1}*\mathcal{A}^{T})$. Then
\begin{align}
\mathop{\max}_{i=1,\cdots,q}f_{i}(\mathcal{X}^{t})\geq \delta_{\infty}^{2}(\mathcal{Q},\boldsymbol{\mathcal{S}})\|\mathcal{X}^{t}-\mathcal{X}^{\star}\|_{F(\mathcal{Q})}^{2},\label{sec4.2e3}\\
\mathbf{E}_{i\sim \mathbf{p}}[f_{i}(\mathcal{X}^{t})]\geq \delta_{\mathbf{p}}^{2}(\mathcal{Q},\boldsymbol{\mathcal{S}})\|\mathcal{X}^{t}-\mathcal{X}^{\star}\|_{F(\mathcal{Q})}^{2}.\label{sec4.2e4}
\end{align}
\end{lemma}
\emph{Proof:} Since $\mathcal{X}^{0}\in \mathbf{Colsp}(\mathcal{Q}^{-1}*\mathcal{A}^{T})$, we have $\mathcal{X}^{t}-\mathcal{X}^{\star}\in \mathbf{Colsp}(\mathcal{Q}^{-1}*\mathcal{A}^{T})$ and consequently
\begin{align*}
\frac{\mathop{\max}\limits_{i=1,\cdots,q}f_{i}(\mathcal{X}^{t})}{\|\mathcal{X}^{t}-\mathcal{X}^{\star}\|_{F(\mathcal{Q})}^{2}}
&=\frac{\mathop{\max}\limits_{i=1,\cdots,q}\|\mathcal{Q}^{\frac{1}{2}}*(\mathcal{X}^{t}-\mathcal{X}^{\star})\|_{F(\mathcal{Z}_{i})}^{2}}{\|\mathcal{X}^{t}-\mathcal{X}^{\star}\|_{F(\mathcal{Q})}^{2}}=\frac{\mathop{\max}\limits_{i=1,\cdots,q}\sum_{j=1}^{p}\|\mathcal{Q}^{\frac{1}{2}}*(\mathcal{X}^{t}-\mathcal{X}^{\star})_{(:,j,:)}\|_{\mathcal{Z}_{i}}^{2}}{\sum_{j=1}^{p}\|(\mathcal{X}^{t}-\mathcal{X}^{\star})_{(:,j,:)}\|_{\mathcal{Q}}^{2}}\\
&\geq\mathop{\min}_{\overrightarrow{\mathcal{V}}\in \mathbf{Range}(\mathcal{Q}^{-1}*\mathcal{A}^{T})}\frac{\mathop{\max}\limits_{i=1,\cdots,q}\sum_{j=1}^{p}\|\mathcal{Q}^{\frac{1}{2}}*\overrightarrow{\mathcal{V}}\|_{\mathcal{Z}_{i}}^{2}}{\sum_{j=1}^{p}\|\overrightarrow{\mathcal{V}}\|_{\mathcal{Q}}^{2}}\\
&\geq\mathop{\min}_{\overrightarrow{\mathcal{V}}\in \mathbf{Range}(\mathcal{Q}^{-1}*\mathcal{A}^{T})}\mathop{\max}_{i=1,\cdots,q}\frac{\|\mathcal{Q}^{\frac{1}{2}}*\overrightarrow{\mathcal{V}}\|_{\mathcal{Z}_{i}}^{2}}{\|\overrightarrow{\mathcal{V}}\|_{\mathcal{Q}}^{2}}=\delta_{\infty}^{2}(\mathcal{Q},\boldsymbol{\mathcal{S}}),~~~\forall t.
\end{align*}
Similarly, we have 
\begin{align*}
\frac{\mathbf{E}_{i\sim \mathbf{p}}[f_{i}(\mathcal{X}^{t})]}{\|\mathcal{X}^{t}-\mathcal{X}^{\star}\|_{F(\mathcal{Q})}^{2}}
&=\frac{\mathbf{E}_{i\sim \mathbf{p}}[\|\mathcal{Q}^{\frac{1}{2}}*(\mathcal{X}^{t}-\mathcal{X}^{\star})\|_{F(\mathcal{Z}_{i})}^{2}]}{\|\mathcal{X}^{t}-\mathcal{X}^{\star}\|_{F(\mathcal{Q})}^{2}}
=\frac{\mathbf{E}_{i\sim \mathbf{p}}[\sum_{j=1}^{p}\|\mathcal{Q}^{\frac{1}{2}}*(\mathcal{X}^{t}-\mathcal{X}^{\star})_{(:,j.:)}\|_{\mathcal{Z}_{i}}^{2}]}{\sum_{j=1}^{p}\|(\mathcal{X}^{t}-\mathcal{X}^{\star})_{(:,j,:)}\|_{\mathcal{Q}}^{2}}\\
&\geq\mathop{\min}_{\overrightarrow{\mathcal{V}}\in \mathbf{Range}(\mathcal{Q}^{-1}*\mathcal{A}^{T})}\frac{\mathbf{E}_{i\sim \mathbf{p}}[\sum_{j=1}^{p}\|\mathcal{Q}^{\frac{1}{2}}*\overrightarrow{\mathcal{V}}\|_{\mathcal{Z}_{i}}^{2}]}{\sum_{j=1}^{p}\|\overrightarrow{\mathcal{V}}\|_{\mathcal{Q}}^{2}}\\
&\geq\mathop{\min}_{\overrightarrow{\mathcal{V}}\in \mathbf{Range}(\mathcal{Q}^{-1}*\mathcal{A}^{T})}\frac{\mathbf{E}_{i\sim \mathbf{p}}[\|\mathcal{Q}^{\frac{1}{2}}*\overrightarrow{\mathcal{V}}\|_{\mathcal{Z}_{i}}^{2}]}{\|\overrightarrow{\mathcal{V}}\|_{\mathcal{Q}}^{2}}=\delta_{\mathbf{p}}^{2}(\mathcal{Q},\boldsymbol{\mathcal{S}}),~~~\forall t.
\end{align*}
Then, the desired results hold.

\begin{lemma}\label{sec4.2lem1}
Let $\mathbf{p}\in \triangle_{q}$ and the set of sketching tubal matrices $\boldsymbol{\mathcal{S}}=[\mathcal{S}_{1},\cdots,\mathcal{S}_{q}]$ be such that $\mathbf{E}_{i\sim \mathbf{p}}[\mathcal{Z}_{i}]$ is T-symmetric T-positive definite with probability $1$. Then 
\begin{align}\label{sec4.2e4-11}
0<\lambda_{\min}(\mathbf{E}_{i\sim \mathbf{p}}[bcirc(\mathcal{Z}_{i})])=\delta_{\mathbf{p}}^{2}(\mathcal{Q},\boldsymbol{\mathcal{S}})\leq\delta_{\infty}^{2}(\mathcal{Q},\boldsymbol{\mathcal{S}})\leq1.
\end{align}
\end{lemma}
\emph{Proof:}
Using $\mathbf{E}_{i\sim \mathbf{p}}[\mathcal{Z}_{i}]=\mathcal{Q}^{-\frac{1}{2}}*\mathcal{A}^{T}*\mathbf{E}_{i\sim \mathbf{p}}[\mathcal{G}_{i}]*\mathcal{A}*\mathcal{Q}^{-\frac{1}{2}}$, as well as the fact that $\mathbf{E}_{i\sim \mathbf{p}}[\mathcal{Z}_{i}]$ is T-symmetric T-positive definite with probability $1$, we obtain 
$$\mathbf{Range}(\mathcal{Q}^{-\frac{1}{2}}*\mathcal{A}^{T})=\mathbf{Range}(\mathbf{E}_{i\sim \mathbf{p}}[\mathcal{Z}_{i}])=\mathbb{K}^{n}_{l}.$$
Hence,
\begin{align*}
\delta_{\mathbf{p}}^{2}(\mathcal{Q},\boldsymbol{\mathcal{S}})&=\mathop{\min}_{\overrightarrow{\mathcal{V}}\in \mathbf{Range}(\mathcal{Q}^{-1}*\mathcal{A}^{T})}\frac{\|\mathcal{Q}^{\frac{1}{2}}*\overrightarrow{\mathcal{V}}\|_{\mathbf{E}_{i\sim \mathbf{p}}[\mathcal{Z}_{i}]}^{2}}{\|\overrightarrow{\mathcal{V}}\|_{\mathcal{Q}}^{2}}
=\mathop{\min}_{\mathcal{Q}^{\frac{1}{2}}*\overrightarrow{\mathcal{V}}\in \mathbf{Range}(\mathbf{E}_{i\sim \mathbf{p}}[\mathcal{Z}_{i}])}\frac{\|\mathcal{Q}^{\frac{1}{2}}*\overrightarrow{\mathcal{V}}\|_{\mathbf{E}_{i\sim \mathbf{p}}[\mathcal{Z}_{i}]}^{2}}{\|\mathcal{Q}^{\frac{1}{2}}*\overrightarrow{\mathcal{V}}\|_{F}^{2}}\\
&=\mathop{\min}_{\text{unfold}(\mathcal{Q}^{\frac{1}{2}}*\overrightarrow{\mathcal{V}})\in \mathbb{R}^{nl}}\frac{\|\text{unfold}(\mathcal{Q}^{\frac{1}{2}}*\overrightarrow{\mathcal{V}})\|_{\text{bcirc}(\mathbf{E}_{i\sim \mathbf{p}}[\mathcal{Z}_{i}])}^{2}}{\|\text{unfold}(\mathcal{Q}^{\frac{1}{2}}*\overrightarrow{\mathcal{V}})\|_{2}^{2}}\\
&=\lambda_{\min}(\text{bcirc}(\mathbf{E}_{i\sim \mathbf{p}}[\mathcal{Z}_{i}]))
=\lambda_{\min}(\mathbf{E}_{i\sim \mathbf{p}}[\text{bcirc}(\mathcal{Z}_{i})])>0,
\end{align*}
and
\begin{align*}
\delta_{\mathbf{p}}^{2}(\mathcal{Q},\boldsymbol{\mathcal{S}})&=\mathop{\min}_{\overrightarrow{\mathcal{V}}\in \mathbf{Range}(\mathcal{Q}^{-1}*\mathcal{A}^{T})}\frac{\|\mathcal{Q}^{\frac{1}{2}}*\overrightarrow{\mathcal{V}}\|_{\mathbf{E}_{i\sim \mathbf{p}}[\mathcal{Z}_{i}]}^{2}}{\|\overrightarrow{\mathcal{V}}\|_{\mathcal{Q}}^{2}}
=\mathop{\min}_{\overrightarrow{\mathcal{V}}\in \mathbf{Range}(\mathcal{Q}^{-1}*\mathcal{A}^{T})}\frac{\mathbf{E}_{i\sim \mathbf{p}}[\|\mathcal{Q}^{\frac{1}{2}}*\overrightarrow{\mathcal{V}}\|_{\mathcal{Z}_{i}}^{2}]}{\|\overrightarrow{\mathcal{V}}\|_{\mathcal{Q}}^{2}}\\
&\leq\mathop{\min}_{\overrightarrow{\mathcal{V}}\in \mathbf{Range}(\mathcal{Q}^{-1}*\mathcal{A}^{T})}\mathop{\max}_{i=1,2,\cdots,q}\frac{\|\mathcal{Q}^{\frac{1}{2}}*\overrightarrow{\mathcal{V}}\|_{\mathcal{Z}_{i}}^{2}}{\|\overrightarrow{\mathcal{V}}\|_{\mathcal{Q}}^{2}}
=\delta_{\infty}^{2}(\mathcal{Q},\boldsymbol{\mathcal{S}}).
\end{align*}
Finally, since the tubal matrix $\mathcal{Z}_{i}$ is an orthogonal projector, we have
\begin{align*}
\delta_{\infty}^{2}(\mathcal{Q},\boldsymbol{\mathcal{S}})&\leq\mathop{\max}_{i=1,2,\cdots,q}\frac{\|\mathcal{Q}^{\frac{1}{2}}*\overrightarrow{\mathcal{V}}\|_{\mathcal{Z}_{i}}^{2}}{\|\overrightarrow{\mathcal{V}}\|_{\mathcal{Q}}^{2}}
=\mathop{\max}_{i=1,2,\cdots,q}\frac{\|\mathcal{Z}_{i}*\mathcal{Q}^{\frac{1}{2}}*\overrightarrow{\mathcal{V}}\|_{F}^{2}}{\|\mathcal{Q}^{\frac{1}{2}}*\overrightarrow{\mathcal{V}}\|_{F}^{2}}
\leq\mathop{\max}_{i=1,2,\cdots,q}\frac{\|\mathcal{Q}^{\frac{1}{2}}*\overrightarrow{\mathcal{V}}\|_{F}^{2}}{\|\mathcal{Q}^{\frac{1}{2}}*\overrightarrow{\mathcal{V}}\|_{F}^{2}}=1.
\end{align*}
Then, the desired results hold.

Next, we give the convergence guarantees of the NTSP, ATSP-MD, ATSP-PR and ATSP-CS methods in turn.

\begin{theorem}\label{sec4.2thm4.1}
Let $\mathcal{X}^{\star}$ satisfy $\mathcal{A}*\mathcal{X}^{\star}=\mathcal{B}$ and $\mathcal{X}^{t}$ be the $t$-th approximation of $\mathcal{X}^{\star}$ calculated by the NTSP method, i.e., Algorithm \ref{NTSP}, with initial iterate $\mathcal{X}^{0}\in \mathbf{Colsp}(\mathcal{Q}^{-1}*\mathcal{A}^{T})$. 
Then $$\mathbf{E}[\|\mathcal{X}^{t}-\mathcal{X}^{\star}\|_{F(\mathcal{Q})}^{2}|\mathcal{X}^{0}]\leq(1-\delta_{\mathbf{p}}^{2}(\mathcal{Q},\boldsymbol{\mathcal{S}}))^{t}\|\mathcal{X}^{0}-\mathcal{X}^{\star}\|_{F(\mathcal{Q})}^{2},$$
where $\delta_{\mathbf{p}}^{2}(\mathcal{Q},\boldsymbol{\mathcal{S}})$ is as defined in (\ref{sec4.2e2}).
\end{theorem}
\emph{Proof:}
From (\ref{sec4.1e2}) and (\ref{sec4.2e4}), we have
\begin{align*}
\mathbf{E}[\|\mathcal{X}^{t+1}-\mathcal{X}^{\star}\|_{F(\mathcal{Q})}^{2}|\mathcal{X}^{t}]
&=\|\mathcal{X}^{t}-\mathcal{X}^{\star}\|_{F(\mathcal{Q})}^{2}-\mathbf{E}_{i^{t}\sim \mathbf{p}}[f_{i^{t}}(\mathcal{X}^{t})]
\leq(1-\delta_{\mathbf{p}}^{2}(\mathcal{Q},\boldsymbol{\mathcal{S}}))\|\mathcal{X}^{t}-\mathcal{X}^{\star}\|_{F(\mathcal{Q})}^{2}.
\end{align*}
Taking the full expectation and unrolling the recurrence, we arrive at this theorem.

\begin{remark}\rm
Since $\delta_{\mathbf{p}}^{2}(\mathcal{Q},\boldsymbol{\mathcal{S}})=\lambda_{\text{min}}(\mathbf{E}_{i\sim \mathbf{p}}[\text{bcirc}(\mathcal{Z}_{i})])$, the conclusion of Theorem \ref{sec4.2thm4.1} can be rewritten as
$$\mathbf{E}[\|\mathcal{X}^{t}-\mathcal{X}^{\star}\|_{F(\mathcal{Q})}^{2}|\mathcal{X}^{0}]\leq(1-\lambda_{\text{min}}(\mathbf{E}_{i\sim \mathbf{p}}[\text{bcirc}(\mathcal{Z}_{i})]))^{t}\|\mathcal{X}^{0}-\mathcal{X}^{\star}\|_{F(\mathcal{Q})}^{2},$$
which is consistent with Theorem \ref{sec3.1thm3.1}.
\end{remark}

\begin{theorem}\label{sec4.2thm4.2}
Let $\mathcal{X}^{\star}$ satisfy $\mathcal{A}*\mathcal{X}^{\star}=\mathcal{B}$ and $\mathcal{X}^{t}$ be the $t$-th approximation of $\mathcal{X}^{\star}$ calculated by the ATSP-MD method, i.e., Algorithm \ref{ATSP-MD}, with initial iterate $\mathcal{X}^{0}\in \mathbf{Colsp}(\mathcal{Q}^{-1}*\mathcal{A}^{T})$. Then $$\|\mathcal{X}^{t}-\mathcal{X}^{\star}\|_{F(\mathcal{Q})}^{2}\leq(1-\delta_{\infty}^{2}(\mathcal{Q},\boldsymbol{\mathcal{S}}))^{t}\|\mathcal{X}^{0}-\mathcal{X}^{\star}\|_{F(\mathcal{Q})}^{2},$$
where $\delta_{\infty}^{2}(\mathcal{Q},\boldsymbol{\mathcal{S}})$ is as defined in (\ref{sec4.2e1}).
\end{theorem}
\emph{Proof:}\emph{Proof:} In view of (\ref{sec4.1e3}) and (\ref{sec4.2e3}), we have
\begin{align*}
\|\mathcal{X}^{t+1}-\mathcal{X}^{\star}\|_{F(\mathcal{Q})}^{2}
&=\|\mathcal{X}^{t}-\mathcal{X}^{\star}\|_{F(\mathcal{Q})}^{2}-\mathop{\max}_{i^{t}=1,2,\dots,q}f_{i^{t}}(\mathcal{X}^{t})
\leq(1-\delta_{\infty}^{2}(\mathcal{Q},\boldsymbol{\mathcal{S}}))\|\mathcal{X}^{t}-\mathcal{X}^{\star}\|_{F(\mathcal{Q})}^{2}.
\end{align*}
Unrolling the recurrence gives this theorem.

\begin{remark}\rm
    Since $\delta_{\mathbf{p}}^{2}(\mathcal{Q},\boldsymbol{\mathcal{S}}) \leq \delta_{\infty}^{2}(\mathcal{Q},\boldsymbol{\mathcal{S}})$, the convergence guarantee for the ATSP-MD method is better than that for the NTSP method.
\end{remark}

\begin{theorem}\label{sec4.2thm4.3}
Let $\mathbf{u}=(\frac{1}{q},\cdots,\frac{1}{q})\in \triangle_{q}$ and $\delta_{\mathbf{u}}^{2}(\mathcal{Q},\boldsymbol{\mathcal{S}})$ be as defined in (\ref{sec4.2e2}).
Let $\mathcal{X}^{\star}$ satisfy $\mathcal{A}*\mathcal{X}^{\star}=\mathcal{B}$ and $\mathcal{X}^{t}$ be the $t$-th approximation of $\mathcal{X}^{\star}$  calculated by the ATSP-PR method, i.e., Algorithm \ref{ATSP-PR}, with initial iterate $\mathcal{X}^{0}\in \mathbf{Colsp}(\mathcal{Q}^{-1}*\mathcal{A}^{T})$. 
Then, for $t\geq 1$, $$\mathbf{E}\left[\|\mathcal{X}^{t+1}-\mathcal{X}^{\star}\|_{F(\mathcal{Q})}^{2}|\mathcal{X}^{t}\right]\leq\left(1-(1+q^{2}\mathbf{Var}_{i\sim \mathbf{u}}[p_{i}^{t}])\delta_{\mathbf{u}}^{2}(\mathcal{Q},\boldsymbol{\mathcal{S}})\right)\|\mathcal{X}^{t}-\mathcal{X}^{\star}\|_{F(\mathcal{Q})}^{2},$$
where $\mathbf{Var}_{i\sim \mathbf{u}}[\cdot]$ denotes the variance taken with respect to the uniform distribution $\mathbf{u}$, i.e.,
$$\mathbf{Var}_{i\sim \mathbf{u}}[v_{i}]=\frac{1}{q}\sum_{i=1}^{q}\left(v_{i}-\frac{1}{q}\sum_{s=1}^{q}v_{s}\right)^{2},~~~\forall ~v\in \mathbb{R}^{q}.$$
Furthermore, 
$$\mathbf{E}\left[\|\mathcal{X}^{t+1}-\mathcal{X}^{\star}\|_{F(\mathcal{Q})}^{2}|\mathcal{X}^{1}\right]\leq\left(1-(1+\frac{1}{q})\delta_{\mathbf{u}}^{2}(\mathcal{Q},\boldsymbol{\mathcal{S}})\right)^{t}\mathbf{E}\left[\|\mathcal{X}^{1}-\mathcal{X}^{\star}\|_{F(\mathcal{Q})}^{2}|\mathcal{X}^{0}\right].$$
\end{theorem}
\emph{Proof:}
First note that, for $i\sim \mathbf{u}$, we have 
\begin{align}\label{sec3e33333}
\mathbf{Var}_{i\sim \mathbf{u}}\left[f_{i}(\mathcal{X}^{t})\right]&=\mathbf{E}_{i\sim \mathbf{u}}\left[f_{i}(\mathcal{X}^{t})^{2}\right]-\mathbf{E}_{i\sim \mathbf{u}}\left[f_{i}(\mathcal{X}^{t})\right]^{2}
=\frac{1}{q}\sum_{i=1}^{q}\left(f_{i}(\mathcal{X}^{t})\right)^{2}-\frac{1}{q^{2}}\left(\sum_{i=1}^{q}f_{i}(\mathcal{X}^{t})\right)^{2}.
\end{align}
Then from (\ref{sec4.2e4}), (\ref{sec3e33333}) and the definition of $\mathbf{p}^{t}$ in Algorithm \ref{ATSP-PR}, we get
\begin{align}
\mathbf{E}_{i\sim \mathbf{p}^{t}}[f_{i}(\mathcal{X}^{t})]
&=\sum_{i=1}^{q}p_{i}^{t}f_{i}(\mathcal{X}^{t}) \notag
=\sum_{i=1}^{q}\frac{(f_{i}(\mathcal{X}^{t}))^{2}}{\sum\limits_{i=1}^{q}f_{i}(\mathcal{X}^{t})}
=\frac{1}{\sum\limits_{i=1}^{q}f_{i}(\mathcal{X}^{t})}\left(q\mathbf{Var}_{i\sim \mathbf{u}}[f_{i}(\mathcal{X}^{t})]+\frac{1}{q}\left(\sum\limits_{i=1}^{q}f_{i}(\mathcal{X}^{t})\right)^{2}\right)\\  \notag
&=\left(1+q^{2}\mathbf{Var}_{i\sim \mathbf{u}}\left[\frac{f_{i}(\mathcal{X}^{t})}{\sum\limits_{i=1}^{q}f_{i}(\mathcal{X}^{t})}\right]\right)\frac{1}{q}\sum\limits_{i=1}^{q}f_{i}(\mathcal{X}^{t})
=(1+q^{2}\mathbf{Var}_{i\sim \mathbf{u}}[p_{i}^{t}])\mathbf{E}_{i\sim \mathbf{u}}[f_{i}(\mathcal{X}^{t})]\\ 
&\geq(1+q^{2}\mathbf{Var}_{i\sim \mathbf{u}}[p_{i}^{t}])\delta_{\mathbf{u}}^{2}(\mathcal{Q},\boldsymbol{\mathcal{S}})\|\mathcal{X}^{t}-\mathcal{X}^{\star}\|_{F(\mathcal{Q})}^{2}.\label{sec4.2e4-1} 
\end{align}
Thus, substituting (\ref{sec4.2e4-1}) into (\ref{sec4.1e2}), we obtain
\begin{align}
\mathbf{E}[\|\mathcal{X}^{t+1}-\mathcal{X}^{\star}\|_{F(\mathcal{Q})}^{2}|\mathcal{X}^{t}]
&=\|\mathcal{X}^{t}-\mathcal{X}^{\star}\|_{F(\mathcal{Q})}^{2}-\mathbf{E}_{i\sim \mathbf{p}^{t}}[f_{i}(\mathcal{X}^{t})]\notag\\
&\leq\left(1-(1+q^{2}\mathbf{Var}_{i\sim \mathbf{\mathbf{u}}}[p_{i}^{t}])\delta_{\mathbf{u}}^{2}(\mathcal{Q},\boldsymbol{\mathcal{S}})\right)\|\mathcal{X}^{t}-\mathcal{X}^{\star}\|_{F(\mathcal{Q})}^{2}.\label{1000}
\end{align}

Next, we further give a lower bound for $\mathbf{Var}_{i\sim \mathbf{u}}[p_{i}^{t}]$. Since
\begin{align*}
\mathcal{Z}_{i^{t}}*\mathcal{Q}^{\frac{1}{2}}*(\mathcal{X}^{t+1}-\mathcal{X}^{\star})
&=\mathcal{Z}_{i^{t}}*\mathcal{Q}^{\frac{1}{2}}*(\mathcal{X}^{t}-\mathcal{Q}^{-\frac{1}{2}}*\mathcal{Z}_{i^{t}}*\mathcal{Q}^{\frac{1}{2}}*(\mathcal{X}^{t}-\mathcal{X}^{\star})-\mathcal{X}^{\star})\\
&=\mathcal{Z}_{i^{t}}*\mathcal{Q}^{\frac{1}{2}}*(\mathcal{X}^{t}-\mathcal{X}^{\star})-\mathcal{Z}_{i^{t}}*\mathcal{Z}_{i^{t}}*\mathcal{Q}^{\frac{1}{2}}*(\mathcal{X}^{t}-\mathcal{X}^{\star})
=\emph{O},
\end{align*}
it follows that
\begin{align}
f_{i^{t}}(\mathcal{X}^{t+1})
&=\|\mathcal{Q}^{\frac{1}{2}}*(\mathcal{X}^{t+1}-\mathcal{X}^{\star})\|_{F(\mathcal{Z}_{i^{t}})}^{2} 
=\|\text{bcirc}(\mathcal{Q}^{\frac{1}{2}})\text{unfold}(\mathcal{X}^{t+1}-\mathcal{X}^{\star})\|_{F(\text{bcirc}(\mathcal{Z}_{i^{t}}))}^{2}\notag\\
&=\sum_{j=1}^{p}\|\text{bcirc}(\mathcal{Q}^{\frac{1}{2}})\text{unfold}(\mathcal{X}^{t+1}-\mathcal{X}^{\star})_{(:,j)}\|_{\text{bcirc}(\mathcal{Z}_{i^{t}})}^{2}\notag\\
&=\sum_{j=1}^{p}\left\langle \text{unfold}(\mathcal{Z}_{i^{t}}*\mathcal{Q}^{\frac{1}{2}}*(\mathcal{X}^{t+1}-\mathcal{X}^{\star}))_{(:,j)},\text{unfold}(\mathcal{Q}^{\frac{1}{2}}*(\mathcal{X}^{t+1}-\mathcal{X}^{\star}))_{(:,j)}\right\rangle
=0,~\forall~t\geq 0,\notag 
\end{align}
which implies $p_{i^{t}}^{t+1}=0$, and hence 
\begin{align}
 \mathbf{Var}_{i\sim \mathbf{u}}[p_{i}^{t+1}]&=\frac{1}{q}\sum_{i=1}^{q}\left(p_{i}^{t+1}-\frac{1}{q}\sum_{s=1}^{q}p_{s}^{t+1}\right)^{2}  
=\frac{1}{q}\sum_{i=1}^{q}\left(p_{i}^{t+1}-\frac{1}{q}\right)^{2}
\geq\frac{1}{q}\left(p_{i^{t}}^{t+1}-\frac{1}{q}\right)^{2}
=\frac{1}{q^{3}}. \label{sec4.2e6}
\end{align}
Therefore, plugging (\ref{sec4.2e6}) into (\ref{1000}), we get
$$\mathbf{E}[\|\mathcal{X}^{t+1}-\mathcal{X}^{\star}\|_{F(\mathcal{Q})}^{2}|\mathcal{X}^{t}]\leq\left(1-(1+\frac{1}{q})\delta_{\mathbf{u}}^{2}(\mathcal{Q},\boldsymbol{\mathcal{S}})\right)\|\mathcal{X}^{t}-\mathcal{X}^{\star}\|_{F(\mathcal{Q})}^{2}.$$
Taking the expectation and unrolling the recursion give this theorem.

\begin{remark}\rm
The convergence rate for the ATSP-PR method is smaller than that for the NTSP method with respect to uniform sampling, and how much smaller depends on the value of $1+q^{2}\mathbf{Var}_{i\sim \mathbf{u}}[p_{i}^{t}]$.
\end{remark}

\begin{theorem}\label{sec4.2thm4.4}
Let $\mathcal{X}^{\star}$ satisfy $\mathcal{A}*\mathcal{X}^{\star}=\mathcal{B}$ and $\mathcal{X}^{t}$ be the $t$-th approximation of $\mathcal{X}^{\star}$ calculated by the ATSP-CS method, i.e., Algorithm \ref{ATSP-CS}, with initial iterate $\mathcal{X}^{0}\in \mathbf{colsp}(\mathcal{Q}^{-1}*\mathcal{A}^{T})$. 
Then
$$\mathbf{E}[\|\mathcal{X}^{t}-\mathcal{X}^{\star}\|_{F(\mathcal{Q})}^{2}|\mathcal{X}^{0}]\leq\left(1-\theta\delta_{\infty}^{2}(\mathcal{Q},\boldsymbol{\mathcal{S}})-(1-\theta)\delta_{\mathbf{p}}^{2}(\mathcal{Q},\boldsymbol{\mathcal{S}})\right)^{t}\|\mathcal{X}^{0}-\mathcal{X}^{\star}\|_{F(\mathcal{Q})}^{2},$$
 where $\delta_{\infty}^{2}(\mathcal{Q},\boldsymbol{\mathcal{S}})$ and $\delta_{\mathbf{p}}^{2}(\mathcal{Q},\boldsymbol{\mathcal{S}})$ are as defined in (\ref{sec4.2e1}) and (\ref{sec4.2e2}), respectively.
\end{theorem}
\emph{Proof:}
Due to $$\mathop{\max}\limits_{j=1,2,\cdots,q}f_{j}(\mathcal{X}^{t})\geq\mathbf{E}_{j\sim \mathbf{p}}[f_{j}(\mathcal{X}^{t})],$$ we know that $\mathfrak{W}_{t}$ defined in (\ref{sec2.e11111}) is not empty 
and $\mathop{\arg\max}\limits_{j=1,2,\cdots,q}f_{j}(\mathcal{X}^{t})\in \mathfrak{W}_{t}$. From the definition of $\mathbf{p}^{t}$ in (\ref{sec2.e22222}), we have $p_{i}^{t}=0$ for all $i\not \in \mathfrak{W}_{k}$, and thus
\begin{align*}
\mathbf{E}[\|\mathcal{X}^{t+1}-\mathcal{X}^{\star}\|_{F(\mathcal{Q})}^{2}|\mathcal{X}^{t}]&=\|\mathcal{X}^{t}-\mathcal{X}^{\star}\|_{F(\mathcal{Q})}^{2}-\mathbf{E}_{i\sim \mathbf{p}^{t}}[f_{i}(\mathcal{X}^{t})]
=\|\mathcal{X}^{t}-\mathcal{X}^{\star}\|_{F(\mathcal{Q})}^{2}-\sum_{i\in \mathfrak{W}_{t}}p_{i}^{t}f_{i}(\mathcal{X}^{t}).
\end{align*}
Note that
\begin{align*}
\sum_{i\in \mathfrak{W}_{t}}f_{i}(\mathcal{X}^{t})p_{i}^{t} \notag
&\geq\sum_{i\in \mathfrak{W}_{t}}\left(\theta\mathop{\max}_{j=1,2,\cdots,q}f_{j}(\mathcal{X}^{t})+(1-\theta)\mathbf{E}_{j\sim \mathbf{p}}[f_{j}(\mathcal{X}^{t})]\right)p_{i}^{t}\\
&=\theta\mathop{\max}_{j=1,2,\cdots,q}f_{j}(\mathcal{X}^{t})+(1-\theta)\mathbf{E}_{j\sim \mathbf{p}}[f_{j}(\mathcal{X}^{t})]\\ 
&\geq(\theta\delta_{\infty}^{2}(\mathcal{Q},\boldsymbol{\mathcal{S}})+(1-\theta)\delta_{\mathbf{p}}^{2}(\mathcal{Q},\boldsymbol{\mathcal{S}}))\|\mathcal{X}^{t}-\mathcal{X}^{\star}\|_{F(\mathcal{Q})}^{2}. 
\end{align*}
Hence, 
\begin{align*}
    \mathbf{E}[\|\mathcal{X}^{t+1}-\mathcal{X}^{\star}\|_{F(\mathcal{Q})}^{2}|\mathcal{X}^{t}]\leq(1-\theta\delta_{\infty}^{2}(\mathcal{Q},\boldsymbol{\mathcal{S}})-(1-\theta)\delta_{\mathbf{p}}^{2}(\mathcal{Q},\boldsymbol{\mathcal{S}}))\|\mathcal{X}^{t}-\mathcal{X}^{\star}\|_{F(\mathcal{Q})}^{2}.
\end{align*}
Taking expectation and unrolling the recursion give this theorem.

\begin{remark}\rm
The convergence rate of the ATSP-CS method is a convex combination of ones of the NTSP and ATSP-MD methods, and hence we can conclude that the closer $\theta$ approaches $0$, the looser the convergence rate of the ATSP-CS method is.
\end{remark}

\begin{remark}\rm
According to Lemma \ref{sec4.2lem1}, we can conclude that the NTSP, ATSP-MD, ATSP-PR and ATSP-CS methods all converge under the assumption that $\mathbf{E}_{i\sim \mathbf{p}}[\mathcal{Z}_{i}]$ is T-symmetric T-positive definite with probability $1$. 
\end{remark}

\begin{remark}\rm
In the TRK setting, using Theorems \ref{sec4.2thm4.1}, \ref{sec4.2thm4.2}, \ref{sec4.2thm4.3}, and \ref{sec4.2thm4.4}, we can get the convergence guarantees for the NTRK, ATRK-MD, ATRK-PR and ATRK-CS methods, respectively.
\end{remark}

\section{The Fourier version of the TSP method}\label{sec:TSPF}
Based on (\ref{sec2e1}) and the discussions following it, we can present an efficient implementation of the TSP method in the Fourier domain, i.e., Algorithm \ref{TSPF}.
\begin{algorithm}[htbp]
\caption{TSP method in Fourier domain}\label{TSPF}
$\textbf{Input:}$ $\mathcal{X}^{0}\in \mathbb{K}^{n\times p}_{l}$, $\mathcal{A}\in \mathbb{K}^{m\times n}_{l}$, $\mathcal{B}\in \mathbb{K}^{m\times p}_{l}$

$\textbf{Parameters:}$ fixed distribution $\mathfrak{D}$ over random tubal matrices, T-symmetric T-positive definite tubal matrix $\mathcal{Q}\in \mathbb{K}^{n\times n}_{l}$

$\widehat{\mathcal{X}}^{0}\longleftarrow \text{fft}(\mathcal{X}^{0},[~],3)$, $\widehat{\mathcal{A}}\longleftarrow
\text{fft}(\mathcal{A},[~],3)$, $\widehat{\mathcal{B}}\longleftarrow \text{fft}(\mathcal{B},[~],3)$, $\widehat{\mathcal{Q}}\longleftarrow \text{fft}(\mathcal{Q},[~],3)$,

$\textbf{for}$ $t=0,1,2,\cdots$

~~~~~Sample an independent copy $\mathcal{S}\sim \mathfrak{D}$

~~~~~$\widehat{\mathcal{S}}\longleftarrow \text{fft}(\mathcal{S},[~],3)$

~~~~~$\textbf{for}$ $k=1,2,\cdots,l$

~~~~~~~~~Compute $\widehat{\mathcal{G}}_{(k)}=\widehat{\mathcal{S}}_{(k)}\left(\widehat{\mathcal{S}}_{(k)}^{H}\widehat{\mathcal{A}}_{(k)}\widehat{\mathcal{Q}}_{(k)}^{-1}\widehat{\mathcal{A}}_{(k)}^{H}\widehat{\mathcal{S}}_{(k)}\right)^{\dag}\widehat{\mathcal{S}}_{(k)}^{H}$

~~~~~~~~~$\widehat{\mathcal{X}}^{t+1}_{(k)}=\widehat{\mathcal{X}}^{t}_{(k)}-{\widehat{\mathcal{Q}}_{(k)}}^{-1}{\widehat{\mathcal{A}}_{(k)}}^{H}\widehat{\mathcal{G}}_{(k)}\left(\widehat{\mathcal{A}}_{(k)}\widehat{\mathcal{X}}_{(k)}^{t}-\widehat{\mathcal{B}}_{(k)}\right)$

~~~~~$\textbf{end for}$

$\textbf{end for}$

$\mathcal{X}^{t+1}\longleftarrow \text{ifft}\left(\widehat{\mathcal{X}}^{t+1},[~],3\right)$

$\textbf{Output:}$ last iterate $\mathcal{X}^{t+1}$
\end{algorithm}

Furthermore, in view of (\ref{sec2e1}), the problem (\ref{sec1e1}) can be reformulated as
\begin{align}\label{sec3.2e1}
\begin{bmatrix}
    \widehat{\mathcal{A}}_{(1)} &   &   &   \\
      & \widehat{\mathcal{A}}_{(2)} &   &   \\
      &   & \ddots &   \\
      &   &   & \widehat{\mathcal{A}}_{(l)} \\
\end{bmatrix}\begin{bmatrix}
                 \widehat{\mathcal{X}}_{(1)} \\
                 \widehat{\mathcal{X}}_{(2)} \\
                 \vdots \\
                 \widehat{\mathcal{X}}_{(l)} \\
\end{bmatrix}=\begin{bmatrix}
                 \widehat{\mathcal{B}}_{(1)} \\
                 \widehat{\mathcal{B}}_{(2)} \\
                 \vdots \\
                 \widehat{\mathcal{B}}_{(l)} \\
\end{bmatrix},
\end{align}
where $\widehat{\mathcal{A}}_{(k)}$, $\widehat{\mathcal{X}}_{(k)}$ and $\widehat{\mathcal{B}}_{(k)}$ for $k=1,2,\cdots,l$ are the frontal slices of $\widehat{\mathcal{A}}=\text{fft}(\mathcal{A},[~],3)$, $\widehat{\mathcal{X}}=\text{fft}(\mathcal{X},[~],3)$ and $\widehat{\mathcal{B}}=\text{fft}(\mathcal{B},[~],3)$, respectively.
As a result, the TSP method 
in the Fourier domain is equivalent to applying the MSP method independently to solve the subsystems $\widehat{\mathcal{A}}_{(k)}\widehat{\mathcal{X}}_{(k)}=\widehat{\mathcal{B}}_{(k)}$ for $k=1,2,\cdots,l$.

We now present a theorem that gives the convergence guarantee for Algorithm \ref{TSPF}.
\begin{theorem}\label{sec3.2thm3.2}
With the notation in Algorithm \ref{TSPF}, 
assume that $\mathbf{E}[\text{bdiag}(\widehat{\mathcal{Z}})]$ is Hermitian positive definite with probability $1$, where $\widehat{\mathcal{Z}}={\rm fft}(\mathcal{Z},[~],3)$ and $\mathcal{Z}$ is as defined in  Theorem \ref{sec3.1thm3.1}. Let $\mathcal{X}^{\star}$ satisfy $\mathcal{A}*\mathcal{X}^{\star}=\mathcal{B}$ and $\mathcal{X}^{t}$ be the $t$-th approximation of $\mathcal{X}^{\star}$ 
with initial iterate $\mathcal{X}^{0}$. Then
\begin{align}\label{sec3.2e2}
\mathbf{E}\left[\|\mathcal{X}^{t}-\mathcal{X}^{\star}\|_{F(\mathcal{Q})}^{2}|\mathcal{X}^{0}\right]\leq\left(1-\mathop{\min}_{k=1,2,\cdots,l}\lambda_{\min}(\mathbf{E}[\widehat{\mathcal{Z}}_{(k)}])\right)^{t}\|\mathcal{X}^{0}-\mathcal{X}^{\star}\|_{F(\mathcal{Q})}^{2}.
\end{align}
\end{theorem}
\emph{Proof:}
According to the properties of t-product, we have the chain of relations
\begin{align*}
\mathbf{E}[\|\mathcal{Z}*\Gamma^{t}\|^{2}_{F}]&=\sum_{j=1}^{p}\mathbf{E}[\langle \text{bcirc}(\mathcal{Z})\text{unfold}(\Gamma^{t})_{(:,j)}, \text{unfold}(\Gamma^{t})_{(:,j)}\rangle]\\
&=\sum_{j=1}^{p}\mathbf{E}[\langle (F_{l}\otimes I_{n})\text{bcirc}(\mathcal{Z})(F^{H}_{l}\otimes I_{n})(F_{l}\otimes I_{n})\text{unfold}(\Gamma^{t})_{(:,j)},(F_{l}\otimes I_{n}) \text{unfold}(\Gamma^{t})_{(:,j)}\rangle]\\
&=\sum_{j=1}^{p}\mathbf{E}[\langle \text{bdiag}(\widehat{\mathcal{Z}})(F_{l}\otimes I_{n})\text{unfold}(\Gamma^{t})_{(:,j)}, (F_{l}\otimes I_{n})\text{unfold}(\Gamma^{t})_{(:,j)}\rangle]\\
&\geq\lambda_{\min}(\mathbf{E}[\text{bdiag}(\widehat{\mathcal{Z}})])\|(F_{l}\otimes I_{n})\text{unfold}(\Gamma^{t})\|_{F}^{2}
=\lambda_{\min}(\mathbf{E}[\text{bdiag}(\widehat{\mathcal{Z}})])\|\Gamma^{t}\|_{F}^{2}\\
&=\mathop{\min}_{k=1,2,\cdots,l}\lambda_{\min}(\mathbf{E}[\widehat{\mathcal{Z}}_{(k)}])\|\Gamma^{t}\|_{F}^{2},
\end{align*}
where the inequality follows from the assumption that $\mathbf{E}[\text{bdiag}(\widehat{\mathcal{Z}})]$ is Hermitian positive definite with probability $1$. Then, we conclude from (\ref{sec3.1e7}) that
\begin{align*}
\mathbf{E}[\|\Gamma^{t+1}\|_{F}^{2}|\mathcal{X}^{t}]
  &=\|\Gamma^{t}\|_{F}^{2}-\mathbf{E}[\|\mathcal{Z}*\Gamma^{t}\|^{2}_{F}]
  \leq\|\Gamma^{t}\|_{F}^{2}-\mathop{\min}_{k=1,2,\cdots,l}\lambda_{\min}(\mathbf{E}[\widehat{\mathcal{Z}}_{(k)}])\|\Gamma^{t}\|_{F}^{2}\\
  &=\left(1-\mathop{\min}_{k=1,2,\cdots,l}\lambda_{\min}(\mathbf{E}[\widehat{\mathcal{Z}}_{(k)}])\right)\|\Gamma^{t}\|_{F}^{2}.
 \end{align*}
That is,
\begin{align*}
  \mathbf{E}\left[\|\mathcal{X}^{t+1}-\mathcal{X}^{\star}\|_{F(\mathcal{Q})}^{2}|\mathcal{X}^{t}\right] \leq \left(1-\mathop{\min}_{k=1,2,\cdots,l}\lambda_{\min}(\mathbf{E}[\widehat{\mathcal{Z}}_{(k)}])\right)\|\mathcal{X}^{t}-\mathcal{X}^{\star}\|_{F(\mathcal{Q})}^{2}.
 \end{align*}
 Taking expectation again and unrolling the recurrence give the result.

\begin{remark}\rm\label{sec4rm4.1}
Similar to Remark \ref{sec3.1rem3.1}, we can obtain that
\begin{align*}
0\leq1-\mathop{\min}_{k=1,2,\cdots,l}\frac{\mathbf{E}[d_{k}]}{n}\leq1-\mathop{\min}_{k=1,2,\cdots,l}\lambda_{\min}(\mathbf{E}[\widehat{\mathcal{Z}}_{(k)}])<1,
\end{align*}
    where $d_{k}=\mathbf{Rank}\left(\widehat{\mathcal{S}}_{(k)}^{H}\widehat{\mathcal{A}}_{(k)}\right)$ for $k=1,2,\cdots,l$. So Algorithm \ref{TSPF} is indeed convergent.
\end{remark}

Next, we give a result in which we consider the random tubal matrix $\mathcal{S}$ with a special discrete probability distribution. To this end, we first recall the definition of the complete discrete sampling matrix presented in \cite{gower2015randomized}: A sampling matrix $S$ is called a complete discrete sampling matrix if it satisfies three conditions, that is, the random matrix $S$ has a discrete distribution, $ S=S_{i}\in \mathbb{C}^{m\times \tau}$ with probability $p_{i}>0$ and make $S_{i}^{H}A$ be of full row rank for $i=1,2,\cdots,q$, and $\boldsymbol{S}=[S_{1},\cdots,S_{q}]\in \mathbb{C}^{m\times q\tau}$ is such that $A^{H}\boldsymbol{S}$ has full row rank.
\begin{corollary}\label{sec3.2cor3.1}
With the notation in Algorithm \ref{TSPF} and Theorem \ref{sec3.2thm3.2}, let $\mathcal{S}$ be a discrete sampling tubal matrix satisfying that $\widehat{\mathcal{S}}_{(k)}$ for $k=1,2,\cdots,l$ are all complete discrete sampling matrices, where $\widehat{\mathcal{S}}={\rm fft
}(\mathcal{S},[~],3)$, and $\mathcal{S}=\mathcal{S}_{i} \in \mathbb{K}^{m}_{l}$ with probability $p_{i}$ for $i=1,2,\cdots,q$. Let $\mathcal{X}^{\star}$ satisfy $\mathcal{A}*\mathcal{X}^{\star}=\mathcal{B}$ and $\mathcal{X}^{t}$ be the $t$-th approximation of $\mathcal{X}^{\star}$ 
with initial iterate $\mathcal{X}^{0}$. Then
when $p_{i}=\frac{\|\mathcal{Q}^{-\frac{1}{2}}*\mathcal{A}^{T}*\mathcal{S}_{i}\|_{F}^{2}}{\|\mathcal{Q}^{-\frac{1}{2}}*\mathcal{A}^{T}*\boldsymbol{\mathcal{S}}\|_{F}^{2}}$ 
with $\boldsymbol{\mathcal{S}}=[\mathcal{S}_{1},\cdots,\mathcal{S}_{q}]$ for $i=1,2,\cdots,q$, we have
\begin{align}\label{sec3.2e4}
\mathbf{E}\left[\left\|\mathcal{X}^{t}-\mathcal{X}^{\star}\right\|_{F(\mathcal{Q})}^{2}|\mathcal{X}^{0}\right]\leq\left(1-\mathop{\min}_{k=1,2,\cdots,l}\frac{\lambda_{\min}\left(\widehat{\boldsymbol{\mathcal{S}}}^{H}_{(k)}\widehat{\mathcal{A}}_{(k)}\widehat{\mathcal{Q}}^{-1}_{(k)}\widehat{\mathcal{A}}^{H}_{(k)}\widehat{\boldsymbol{\mathcal{S}}}_{(k)}\right)}{\|\mathcal{Q}^{-\frac{1}{2}}*\mathcal{A}^{T}*\boldsymbol{\mathcal{S}}\|_{F}^{2}}\right)^{t}\left\|\mathcal{X}^{0}-\mathcal{X}^{\star}\right\|_{F(\mathcal{Q})}^{2};
\end{align}
when $p_{i}=\frac{1}{q}$ for $i=1,2,\cdots,q$, we have
\begin{align}\label{sec3.2e5}
\mathbf{E}\left[\left\|\mathcal{X}^{t}-\mathcal{X}^{\star}\right\|_{F(\mathcal{Q})}^{2}|\mathcal{X}^{0}\right]\leq\left(1-\mathop{\min}_{k=1,2,\cdots,l}\frac{\lambda_{\min}\left(\widehat{\boldsymbol{\mathcal{S}}}_{(k)}^{H}\widehat{\mathcal{A}}_{(k)}\widehat{\mathcal{Q}}^{-1}_{(k)}\widehat{\mathcal{A}}_{(k)}^{H}\widehat{\boldsymbol{\mathcal{S}}}_{(k)}\right)}{q\max\limits_{i=1,2,\cdots,q}\left\|\widehat{\mathcal{Q}}_{(k)}^{-\frac{1}{2}}\widehat{\mathcal{A}}_{(k)}^{H}(\widehat{\mathcal{S}}_{i})_{(k)}\right\|_{F}^{2}}\right)^{t}\left\|\mathcal{X}^{0}-\mathcal{X}^{\star}\right\|_{F(\mathcal{Q})}^{2}.
\end{align}
\end{corollary}
\emph{Proof:}
 Since $\mathcal{S}$ satisfies that $\widehat{\mathcal{S}}_{(k)}$ for $k=1,2,\cdots,l$ are all complete discrete sampling matrices, 
we can get that $\mathbf{E}[\text{bdiag}(\widehat{\mathcal{Z}})]$ is Hermitian positive definite, which implies that such sketching tubal matrix satisfies the assumptions in Theorem \ref{sec3.2thm3.2}. Let $$D_{k}=\text{diag}\left(\sqrt{p_{1}}\left((\widehat{\mathcal{S}}_{1})_{(k)}^{H}\widehat{\mathcal{A}}_{(k)}\widehat{\mathcal{Q}}^{-1}_{(k)}\widehat{\mathcal{A}}_{(k)}^{H}(\widehat{\mathcal{S}}_{1})_{(k)}\right)^{-\frac{1}{2}},\cdots,\sqrt{p_{q}}\left((\widehat{\mathcal{S}}_{q})_{(k)}^{H}\widehat{\mathcal{A}}_{(k)}\widehat{\mathcal{Q}}^{-1}_{(k)}\widehat{\mathcal{A}}_{(k)}^{H}(\widehat{\mathcal{S}}_{q})_{(k)}\right)^{-\frac{1}{2}}\right).$$
Then, $\mathbf{E}[\widehat{\mathcal{Z}}_{(k)}]$ can be expressed as
\begin{align*}
\mathbf{E}[\widehat{\mathcal{Z}}_{(k)}]&=\sum_{i=1}^{q}\left(\widehat{\mathcal{Q}}^{-\frac{1}{2}}_{(k)}\widehat{\mathcal{A}}_{(k)}^{H}(\widehat{\mathcal{S}_{i}})_{(k)}\left((\widehat{\mathcal{S}}_{i})_{(k)}^{H}\widehat{\mathcal{A}}_{(k)}\widehat{\mathcal{Q}}^{-1}_{(k)}\widehat{\mathcal{A}}_{(k)}^{H}(\widehat{\mathcal{S}_{i}})_{(k)}\right)^{-1}(\widehat{\mathcal{S}}_{i})_{(k)}^{H}\widehat{\mathcal{A}}_{(k)}\widehat{\mathcal{Q}}^{-\frac{1}{2}}_{(k)}\right)p_{i}\\
&=\widehat{\mathcal{Q}}^{-\frac{1}{2}}_{(k)}\widehat{\mathcal{A}}_{(k)}^{H}\widehat{\boldsymbol{\mathcal{S}}}_{(k)}D_{k}^{2}\widehat{\boldsymbol{\mathcal{S}}}_{(k)}^{H}\widehat{\mathcal{A}}_{(k)}\widehat{\mathcal{Q}}^{-\frac{1}{2}}_{(k)}.
\end{align*}
Therefore, we obtain
\begin{align*}
\lambda_{\min}\left(\mathbf{E}[\widehat{\mathcal{Z}}_{(k)}]\right)&=\lambda_{\min}\left(\widehat{\mathcal{Q}}^{-\frac{1}{2}}_{(k)}\widehat{\mathcal{A}}_{(k)}^{H}\widehat{\boldsymbol{\mathcal{S}}}_{(k)}D_{k}^{2}\widehat{\boldsymbol{\mathcal{S}}}_{(k)}^{H}\widehat{\mathcal{A}}_{(k)}\widehat{\mathcal{Q}}^{-\frac{1}{2}}_{(k)}\right)
=\lambda_{\min}\left(\widehat{\boldsymbol{\mathcal{S}}}_{(k)}^{H}\widehat{\mathcal{A}}_{(k)}\widehat{\mathcal{Q}}^{-1}_{(k)}\widehat{\mathcal{A}}_{(k)}^{H}\widehat{\boldsymbol{\mathcal{S}}}_{(k)}D_{k}^{2}\right)\\
&\geq\lambda_{\min}\left(\widehat{\boldsymbol{\mathcal{S}}}_{(k)}^{H}\widehat{\mathcal{A}}_{(k)}\widehat{\mathcal{Q}}^{-1}_{(k)}\widehat{\mathcal{A}}_{(k)}^{H}\widehat{\boldsymbol{\mathcal{S}}}_{(k)}\right)\lambda_{\min}\left(D_{k}^{2}\right).
\end{align*}
When $p_{i}=\frac{\|\mathcal{Q}^{-\frac{1}{2}}*\mathcal{A}^{T}*\mathcal{S}_{i}\|_{F}^{2}}{\|\mathcal{Q}^{-\frac{1}{2}}*\mathcal{A}^{T}*\boldsymbol{\mathcal{S}}\|_{F}^{2}}$ with $\boldsymbol{\mathcal{S}}=[\mathcal{S}_{1},\cdots,\mathcal{S}_{q}]$ for $i=1,2,\cdots,q$,
according to the properties of t-product, we have the chain of relations
\begin{align*}
\left\|\mathcal{Q}^{-\frac{1}{2}}*\mathcal{A}^{T}*\mathcal{S}_{i}\right\|_{F}^{2}&=\frac{1}{l}\sum_{k=1}^{l}\left\|\widehat{\mathcal{Q}}^{-\frac{1}{2}}_{(k)}\widehat{\mathcal{A}}_{(k)}^{H}(\widehat{\mathcal{S}}_{i})_{(k)}\right\|_{F}^{2}
=\frac{1}{l}\sum_{k=1}^{l}\mathbf{Tr}\left((\widehat{\mathcal{S}}_{i})_{(k)}^{H}\widehat{\mathcal{A}}_{(k)}\widehat{\mathcal{Q}}^{-1}_{(k)}\widehat{\mathcal{A}}_{(k)}^{H}(\widehat{\mathcal{S}_{i}})_{(k)}\right)\\
&\geq\frac{1}{l}\sum_{k=1}^{l}\lambda_{\max}\left((\widehat{\mathcal{S}}_{i})_{(k)}^{H}\widehat{\mathcal{A}}_{(k)}\widehat{\mathcal{Q}}^{-1}_{(k)}\widehat{\mathcal{A}}_{(k)}^{H}(\widehat{\mathcal{S}_{i}})_{(k)}\right)
=\lambda_{\max}\left((\widehat{\mathcal{S}}_{i})_{(k)}^{H}\widehat{\mathcal{A}}_{(k)}\widehat{\mathcal{Q}}^{-1}_{(k)}\widehat{\mathcal{A}}_{(k)}^{H}(\widehat{\mathcal{S}_{i}})_{(k)}\right),
\end{align*}
which immediately yields
\begin{align*}
\lambda_{\min}(D_{k}^{2})&=\min_{i=1,2,\cdots,q}\left(\frac{\|\mathcal{Q}^{-\frac{1}{2}}*\mathcal{A}^{T}*\mathcal{S}_{i}\|_{F}^{2}}{\|\mathcal{Q}^{-\frac{1}{2}}*\mathcal{A}^{T}*\boldsymbol{\mathcal{S}}\|_{F}^{2}}\cdot\frac{1}{\lambda_{\max}\left((\widehat{\mathcal{S}}_{i})_{(k)}^{H}\widehat{\mathcal{A}}_{(k)}\widehat{\mathcal{Q}}^{-1}_{(k)}\widehat{\mathcal{A}}_{(k)}^{H}(\widehat{\mathcal{S}_{i}})_{(k)}\right)}\right)\\
&\geq\min_{i=1,2,\cdots,q}\frac{1}{\left\|\mathcal{Q}^{-\frac{1}{2}}*\mathcal{A}^{T}*\boldsymbol{\mathcal{S}}\right\|_{F}^{2}}=\frac{1}{\left\|\mathcal{Q}^{-\frac{1}{2}}*\mathcal{A}^{T}*\boldsymbol{\mathcal{S}}\right\|_{F}^{2}}.
\end{align*}
As a consequence,
\begin{align}\label{sec3.2e6}
\lambda_{\min}\left(\mathbf{E}[\widehat{\mathcal{Z}}_{(k)}]\right)
\geq\frac{\lambda_{\min}\left(\widehat{\boldsymbol{\mathcal{S}}}_{(k)}^{H}\widehat{\mathcal{A}}_{(k)}\widehat{\mathcal{Q}}^{-1}_{(k)}\widehat{\mathcal{A}}_{(k)}^{H}\widehat{\boldsymbol{\mathcal{S}}}_{(k)}\right)}{\left\|\mathcal{Q}^{-\frac{1}{2}}*\mathcal{A}^{T}*\boldsymbol{\mathcal{S}}\right\|_{F}^{2}}.
\end{align}
When $p_{i}=\frac{1}{q}$ for $i=1,2,\cdots,q$, we have
\begin{align*}
\lambda_{\min}\left(D_{k}^{2}\right)&=\frac{1}{q}\left(\frac{1}{\max\limits_{i=1,2,\cdots,q}\lambda_{\max}\left((\widehat{\mathcal{S}}_{i})_{(k)}^{H}\widehat{\mathcal{A}}_{(k)}\widehat{\mathcal{Q}}^{-1}_{(k)}\widehat{\mathcal{A}}_{(k)}^{H}(\widehat{\mathcal{S}}_{i})_{(k)}\right)}\right)
\geq\frac{1}{q\max\limits_{i=1,2,\cdots,q}\left\|\widehat{\mathcal{Q}}_{(k)}^{-\frac{1}{2}}\widehat{\mathcal{A}}_{(k)}^{H}(\widehat{\mathcal{S}}_{i})_{(k)}\right\|_{F}^{2}}.
\end{align*}
Thus
\begin{align}\label{sec3.2e7}
\lambda_{\min}\left(\mathbf{E}\left[\widehat{\mathcal{Z}}_{(k)}\right]\right)
\geq\frac{\lambda_{\min}\left(\widehat{\boldsymbol{\mathcal{S}}}_{(k)}^{H}\widehat{\mathcal{A}}_{(k)}\widehat{\mathcal{Q}}^{-1}_{(k)}\widehat{\mathcal{A}}_{(k)}^{H}\widehat{\boldsymbol{\mathcal{S}}}_{(k)}\right)}{q\max\limits_{i=1,2,...,q}\left\|\widehat{\mathcal{Q}}_{(k)}^{-\frac{1}{2}}\widehat{\mathcal{A}}_{(k)}^{H}(\widehat{\mathcal{S}}_{i})_{(k)}\right\|_{F}^{2}}.
\end{align}
Combine (\ref{sec3.2e2}), (\ref{sec3.2e6}) and (\ref{sec3.2e7}) to reach the main results (\ref{sec3.2e4}) and (\ref{sec3.2e5}).

\begin{remark}\rm
In Corollary \ref{sec3.2cor3.1},  
choosing $\mathcal{S}_{i}=\mathcal{I}_{(:,i,:)}\in \mathbb{K}^{m}_{l}$ for $i=1,2,\cdots,m$ and $\mathcal{Q}=\mathcal{I} \in \mathbb{K}^{n\times n }_{l}$, and setting the probability $p_{i}=\frac{\|\mathcal{A}_{(i,:,:)}\|_{F}^{2}}{\|\mathcal{A}\|_{F}^{2}}$ (proportional to the magnitude of $i$-th horizontal slice of $\mathcal{A}$) and $p_{i}=\frac{1}{m}$ (uniform sampling) lead to
\begin{align}\label{sec5e2}
\mathbf{E}\left[\left\|\mathcal{X}^{t}-\mathcal{X}^{\star}\right\|_{F}^{2}|\mathcal{X}^{0}\right]\leq\left(1-\mathop{\min}_{k=1,2,\cdots,l}\frac{\lambda_{\min}\left(\widehat{\mathcal{A}}_{(k)}\widehat{\mathcal{A}}^{H}_{(k)}\right)}{\|\mathcal{A}\|_{F}^{2}}\right)^{t}\left\|\mathcal{X}^{0}-\mathcal{X}^{\star}\right\|_{F}^{2}
\end{align}
and
\begin{align}\label{sec5e1}
\mathbf{E}\left[\left\|\mathcal{X}^{t}-\mathcal{X}^{\star}\right\|_{F}^{2}|\mathcal{X}^{0}\right]\leq\left(1-\mathop{\min}_{k=1,2,\cdots,l}\frac{\lambda_{\min}\left(\widehat{\mathcal{A}}_{(k)}\widehat{\mathcal{A}}_{(k)}^{H}\right)}{m\max\limits_{i=1,2,\cdots,q}\left(\left(\widehat{\mathcal{A}_{(i,:,:)}}\right)_{(k)}\left(\widehat{\mathcal{A}_{(i,:,:)}^{T}}\right)_{(k)}\right)}\right)^{t}\left\|\mathcal{X}^{0}-\mathcal{X}^{\star}\right\|_{F}^{2},
\end{align}
respectively, where 
(\ref{sec5e1}) is just the result for the TRK method 
given in Theorem $4.1$ in \cite{ma2021randomized}.  
\end{remark}
\begin{remark}\rm
For the NTSP and three adaptive TSP methods discussed in Subsection \ref{sec:ATSP}, we can also implement them in the Fourier domain, 
and obtain the corresponding convergence guarantees 
in a similar way. The details are omitted here.
\end{remark}

\section{Two improved strategies}  
\label{sec:TSPsub3}
The sketching tubal matrix $\mathcal{S}$ appearing in the algorithms proposed in Sections \ref{sec:TSP} and \ref{sec:TSPF} can be formed as done in 
\cite{tarzanagh2018fast, ma2021randomized, zhang2018randomizedddd, qiliqun2021tttt}. However, as explained in Section \ref{secint}, in this case,  $\widehat{\mathcal{S}}_{(k)}$, for $k=1,2,\cdots,l$, will be the same, and hence the sketching matrices for all the subsystems $\widehat{\mathcal{A}}_{(k)}\widehat{\mathcal{X}}_{(k)}=\widehat{\mathcal{B}}_{(k)}$, for $k=1,2,\cdots,l$, are the same. 
For complex-valued problems, Ma and Molitor \cite{ma2021randomized} proposed to select different sketching matrices, i.e., select different indices, for different subsystems. However, this strategy doesn't work for real-valued problems considered in this paper. This is because, in this case, 
the approximate solution $\mathcal{X}^{t+1}=\text{ifft}(\widehat{\mathcal{X}}^{t+1},[~],3)$ is no longer real-valued. 
To tackle this problem, 
we propose two improved strategies. 
The first one is based on the following equivalence transformation:
$$\{\mathcal{X}\in \mathbb{R}^{n\times p\times l}~|~\mathcal{A}*\mathcal{X}=\mathcal{B}\}=\left\{\mathcal{X}\in \mathbb{R}^{n\times p\times l}~|~\begin{bmatrix}
                                                                                                                                        \text{Re}(\mathcal{A}) \\
                                                                                                                                        \text{Im}(\mathcal{A}) \\
                                                                                                                                      \end{bmatrix}
*\mathcal{X}=\begin{bmatrix}
                 \text{Re}(\mathcal{B}) \\
                 \text{Im}(\mathcal{B}) \\
          \end{bmatrix}
\right\}.$$
Putting this equivalence transformation into the TSP method, we can get the first improved algorithm, i.e. Algorithm \ref{TSP-1}. And the convergence of the TSP-\uppercase\expandafter{\romannumeral1} method is provided in 
Theorem \ref{sec3.3thm3.3}.

\begin{algorithm}[htbp]
\caption{TSP-\uppercase\expandafter{\romannumeral1} method}\label{TSP-1}
$\textbf{Input:}$ $\mathcal{X}^{0}\in \mathbb{K}^{n\times p}_{l}$, $\mathcal{A}\in \mathbb{K}^{m\times n}_{l}$, $\mathcal{B}\in \mathbb{K}^{m\times p}_{l}$

$\textbf{Parameters:}$ fixed distribution $\mathfrak{D}_{k}$ over random matrices for $k=1,2,\cdots,l$, T-symmetric T-positive definite tubal matrix $\mathcal{Q}\in \mathbb{K}^{n\times n}_{l}$

$\widehat{\mathcal{X}}^{0}\longleftarrow \text{fft}(\mathcal{X}^{0},[~],3)$, $\widehat{\mathcal{A}}\longleftarrow \text{fft}(\mathcal{A},[~],3)$, $\widehat{\mathcal{B}}\longleftarrow \text{fft}(\mathcal{B},[~],3)$, $\widehat{\mathcal{Q}}\longleftarrow \text{fft}(\mathcal{Q},[~],3)$

$\textbf{for}$ $t=0,1,2,\cdots$

~~~~~$\textbf{for}$ $k=1,2,\cdots,l$

~~~~~~~~~$S_{k}\sim \mathfrak{D}_{k}$

~~~~~~~~~$\mathcal{S}_{(k)}=S_{k}$

~~~~~~~~~$(\mathcal{\Check{A}}_\mathcal{S})_{(k)}=\mathcal{S}_{(k)}^{H}\widehat{\mathcal{A}}_{(k)}$, $(\mathcal{\Check{B}}_{\mathcal{S}})_{(k)}=\mathcal{S}_{(k)}^{H}\widehat{\mathcal{B}}_{(k)}$

~~~~~$\textbf{end for}$

~~~~~$\widetilde{\mathcal{A}_{\mathcal{S}}}=\text{ifft}(\mathcal{\Check{A}}_\mathcal{S},[~],3)$, $\widetilde{\mathcal{B}_{\mathcal{S}}}=\text{ifft}(\mathcal{\Check{B}}_\mathcal{S},[~],3)$

~~~~~$\mathcal{A}_{\mathcal{S}}=\begin{bmatrix}
                                   \text{Re}(\widetilde{\mathcal{A}_{\mathcal{S}}}) \\
                                   \text{Im}(\widetilde{\mathcal{A}_{\mathcal{S}}}) \\
                                 \end{bmatrix}
$, $\mathcal{B}_{\mathcal{S}}= \begin{bmatrix}
                                  \text{Re}(\widetilde{\mathcal{B}_{\mathcal{S}}})\\
                                  \text{Im}(\widetilde{\mathcal{B}_{\mathcal{S}}}) \\
                                \end{bmatrix}$

~~~~~$\widehat{\mathcal{A}_{\mathcal{S}}}=\text{fft}(\mathcal{A}_{\mathcal{S}},[~],3)$, $\widehat{\mathcal{B}_{\mathcal{S}}}=\text{fft}(\mathcal{B}_{\mathcal{S}},[~],3)$

~~~~~$\textbf{for}$ $k=1,2,\cdots,l$

~~~~~~~~~$\widehat{\mathcal{X}}^{t+1}_{(k)}=\widehat{\mathcal{X}}^{t}_{(k)}-{\widehat{\mathcal{Q}}_{(k)}}^{-1}(\widehat{\mathcal{A}_{\mathcal{S}}})_{(k)}^{H}\left((\widehat{\mathcal{A}_{\mathcal{S}}})_{(k)}{\widehat{\mathcal{Q}}_{(k)}}^{-1}(\widehat{\mathcal{A}_{\mathcal{S}}})_{(k)}^{H}\right)^{\dag}\left((\widehat{\mathcal{A}_{\mathcal{S}}})_{(k)}\widehat{\mathcal{X}}^{t}_{(k)}-{(\widehat{\mathcal{B}_{\mathcal{S}}})_{(k)}}\right)$

~~~~~$\textbf{end for}$

$\textbf{end for}$

$\mathcal{X}^{t+1}\longleftarrow \text{ifft}(\widehat{\mathcal{X}}^{t+1},[~],3)$

 $\textbf{Output:}$ last iterate $\mathcal{X}^{t+1}$
\end{algorithm}

\begin{theorem}\label{sec3.3thm3.3}
With the notation in Algorithm \ref{TSP-1}, 
assume that $\mathbf{E}[\widehat{\mathcal{Z}}_{(k)}]$ is Hermitian positive definite with probability $1$, where $\widehat{\mathcal{Z}}_{(k)}=\widehat{\mathcal{Q}}^{-\frac{1}{2}}_{(k)}\widehat{\mathcal{A}}^{H}_{(k)}S_{k}(S_{k}^{H}\widehat{\mathcal{A}}_{(k)}\widehat{\mathcal{Q}}^{-1}_{(k)}\widehat{\mathcal{A}}_{(k)}^{H}S_{k})^{-1}S_{k}^{H}\widehat{\mathcal{A}}_{(k)}\widehat{\mathcal{Q}}^{-\frac{1}{2}}_{(k)}$, for $k=1,2,\cdots,l$. Let $\mathcal{X}^{\star}$ satisfy $\mathcal{A}*\mathcal{X}^{\star}=\mathcal{B}$ and $\mathcal{X}^{t}$ be the $t$-th approximation of $\mathcal{X}^{\star}$
with initial iterate $\mathcal{X}^{0}$. Then
$$\mathbf{E}\left[\left\|\mathcal{X}^{t}-\mathcal{X}^{\star}\right\|_{F(\mathcal{Q})}^{2}|\mathcal{X}^{0}\right]\leq\left(1-\mathop{\min}_{k=1,2,\cdots,l}\lambda_{\min}\left(\mathbf{E}[\widehat{\mathcal{Z}}_{(k)}]\right)\right)^{t}\left\|\mathcal{X}^{0}-\mathcal{X}^{\star}\right\|_{F(\mathcal{Q})}^{2}.$$
\end{theorem}
\emph{Proof:}
 The update of Algorithm \ref{TSP-1} can be expressed by
\begin{align*}
\mathcal{X}^{t+1}=\mathcal{X}^{t}-\mathcal{Q}^{-1}*\underline{\mathcal{W}}*(\mathcal{X}^{t}-\mathcal{X}^{\star}),
\end{align*}
where $$\underline{\mathcal{W}}=\begin{bmatrix}
                         \text{Re}(\mathcal{S}^{H}*\mathcal{A}) \\
                         \text{Im}(\mathcal{S}^{H}*\mathcal{A}) \\
\end{bmatrix}^{T}*\left(\begin{bmatrix}
                         \text{Re}(\mathcal{S}^{H}*\mathcal{A}) \\
                         \text{Im}(\mathcal{S}^{H}*\mathcal{A}) \\
                       \end{bmatrix}*\mathcal{Q}^{-1}*\begin{bmatrix}
                         \text{Re}(\mathcal{S}^{H}*\mathcal{A}) \\
                         \text{Im}(\mathcal{S}^{H}*\mathcal{A}) \\
                      \end{bmatrix}^{T}\right)^{\dag}*\begin{bmatrix}
                         \text{Re}(\mathcal{S}^{H}*\mathcal{A}) \\
                         \text{Im}(\mathcal{S}^{H}*\mathcal{A}) \\
                      \end{bmatrix}$$ 
and $\mathcal{S}=\text{ifft}(\widehat{\mathcal{S}},[~],3)$ with $\widehat{\mathcal{S}}$ being a random tubal matrix whose frontal slices are $S_{k}$ for $k=1,2,\cdots,l$.
For $\mathcal{Z}=\text{ifft}(\widehat{\mathcal{Z}},[~],3)$ and $\underline{\mathcal{Z}}=\mathcal{Q}^{-\frac{1}{2}}*\underline{\mathcal{W}}*\mathcal{Q}^{-\frac{1}{2}}$,
it is clear that $\text{bcirc}(\mathcal{Z})$ and $\text{bcirc}(\underline{\mathcal{Z}})$ are both orthogonal projections, and hence the spectrums of $\text{bcirc}(\mathcal{Z})$ and $\text{bcirc}(\underline{\mathcal{Z}})$ are contained in $\{0,1\}$. According to $$\text{bcirc}(\mathcal{S}^{H}*\mathcal{A})=\text{Re}(\text{bcirc}(\mathcal{S}^{H}*\mathcal{A}))+\text{Im}(\text{bcirc}(\mathcal{S}^{H}*\mathcal{A}))i=\begin{bmatrix}
                                                                                                                             I & iI \\
                                                                                                                          \end{bmatrix}\begin{bmatrix}
                                                                                                                                    \text{Re}(\text{bcirc}(\mathcal{S}^{H}*\mathcal{A})) \\
                                                                                                                                    \text{Im}(\text{bcirc}(\mathcal{S}^{H}*\mathcal{A})) \\
                                                                                                                                  \end{bmatrix},
$$
we can conclude that
$$\mathbf{Rank}(\text{bcirc}(\mathcal{S}^{H}*\mathcal{A}))\leq \mathbf{Rank}\left(\begin{bmatrix}
                                                                                                                                    \text{Re}(\text{bcirc}(\mathcal{S}^{H}*\mathcal{A})) \\
                                                                                                                                    \text{Im}(\text{bcirc}(\mathcal{S}^{H}*\mathcal{A})) \\
                                                                                                                                  \end{bmatrix}\right),$$
which implies that 
$\lambda_{\min}(\text{bcirc}(\underline{\mathcal{Z}}))\geq \lambda_{\min}(\text{bcirc}(\mathcal{Z}))\geq 0$. Hence, $\lambda_{\min}(\text{bdiag}(\widehat{\underline{\mathcal{Z}}}))=\lambda_{\min}((F_{l}\otimes I_{n})\text{bcirc}(\underline{\mathcal{Z}})(F_{l}^{H}\otimes I_{n}))\geq \lambda_{\min}((F_{l}\otimes I_{n})\text{bcirc}(\mathcal{Z})(F_{l}^{H}\otimes I_{n}))=\lambda_{\min}(\text{bdiag}(\widehat{\mathcal{Z}}))\geq 0$.
In addition, the hypothesis that $\mathbf{E}[\widehat{\mathcal{Z}}_{(k)}]$ is Hermitian positive definite for $k=1,2,\cdots,l$ implies that $\mathbf{E}[\text{bdiag}(\widehat{\mathcal{Z}})]$ is also Hermitian positive definite. Therefore, $\mathbf{E}[\text{bdiag}(\widehat{\underline{\mathcal{Z}}})]\geq\mathbf{E}[\text{bdiag}(\widehat{\mathcal{Z}})]>O$, which means that $\lambda_{\min}(\mathbf{E}[\text{bdiag}(\widehat{\underline{\mathcal{Z}}})])\geq\lambda_{\min}(\mathbf{E}[\text{bdiag}(\widehat{\mathcal{Z}})])>0$. Thus,
similar to the proof of Theorem \ref{sec3.2thm3.2}, we obtain
\begin{align*}
\mathbf{E}\left[\left\|\mathcal{X}^{t}-\mathcal{X}^{\star}\right\|_{F(\mathcal{Q})}^{2}|\mathcal{X}^{0}\right]&\leq\left(1-\lambda_{\min}(\mathbf{E}[\text{bdiag}(\widehat{\underline{\mathcal{Z}}})])\right)^{t}\left\|\mathcal{X}^{0}-\mathcal{X}^{\star}\right\|_{F(\mathcal{Q})}^{2}\\
&\leq\left(1-\lambda_{\min}(\mathbf{E}[\text{bdiag}(\widehat{\mathcal{Z}})])\right)^{t}\left\|\mathcal{X}^{0}-\mathcal{X}^{\star}\right\|_{F(\mathcal{Q})}^{2}\\
&=\left(1-\mathop{\min}_{k=1,2,\cdots,l}\lambda_{\min}(\mathbf{E}[\widehat{\mathcal{Z}}_{(k)}])\right)^{t}\left\|\mathcal{X}^{0}-\mathcal{X}^{\star}\right\|_{F(\mathcal{Q})}^{2}.
\end{align*}

\begin{remark}\rm
According to Remark \ref{sec4rm4.1}, we can conclude that Algorithm \ref{TSP-1} is convergent.
\end{remark}

\begin{corollary}\label{sec3.3cor3.2}
With the notation in Algorithm \ref{TSP-1} and Theorem \ref{sec3.3thm3.3}, let $S_{k}$ be a complete discrete sampling matrix for $k=1,2,\cdots,l$ and $S_{k}=S_{k_{i}} \in \mathbb{C}^{m}$ with probability $p_{k_{i}}$ for $i=1,2,\cdots,q$. Let $\mathcal{X}^{\star}$ satisfy $\mathcal{A}*\mathcal{X}^{\star}=\mathcal{B}$ and $\mathcal{X}^{t}$ be the $t$-th approximation of $\mathcal{X}^{\star}$
with initial iterate $\mathcal{X}^{0}$. Then when $p_{k_{i}}=\frac{\|\widehat{\mathcal{Q}}_{(k)}^{-\frac{1}{2}}\widehat{\mathcal{A}}_{(k)}^{H}S_{k_{i}}\|_{F}^{2}}{\|\widehat{\mathcal{Q}}_{(k)}^{-\frac{1}{2}}\widehat{\mathcal{A}}^{H}_{(k)}\boldsymbol{S}_{k}\|_{F}^{2}}$ where $\boldsymbol{S}_{k}=[S_{k_{1}},\cdots,S_{k_{q}}]$, we have
\begin{align*}
\mathbf{E}\left[\left\|\mathcal{X}^{t}-\mathcal{X}^{\star}\right\|_{F(\mathcal{Q})}^{2}|\mathcal{X}^{0}\right]\leq\left(1-\mathop{\min}_{k=1,2,\cdots,l}\frac{\lambda_{\min}\left(\boldsymbol{S}_{k}^{H}\widehat{\mathcal{A}}_{(k)}\widehat{\mathcal{Q}}^{-1}_{(k)}\widehat{\mathcal{A}}_{(k)}^{H}\boldsymbol{S}_{k}\right)}{\|\widehat{\mathcal{Q}}_{(k)}^{-\frac{1}{2}}\widehat{\mathcal{A}}^{H}_{(k)}\boldsymbol{S}_{k}\|_{F}^{2}}\right)^{t}\|\mathcal{X}^{0}-\mathcal{X}^{\star}\|_{F(\mathcal{Q})}^{2};
\end{align*}
when $p_{k_{i}}=\frac{1}{q}$, we have
\begin{align*}
\mathbf{E}\left[\left\|\mathcal{X}^{t}-\mathcal{X}^{\star}\right\|_{F(\mathcal{Q})}^{2}|\mathcal{X}^{0}\right]\leq\left(1-\mathop{\min}_{k=1,2,\cdots,l}\frac{\lambda_{\min}\left(\boldsymbol{S}_{k}^{H}\widehat{\mathcal{A}}_{(k)}\widehat{\mathcal{Q}}^{-1}_{(k)}\widehat{\mathcal{A}}_{(k)}^{H}\boldsymbol{S}_{k}\right)}{q\max\limits_{i=1,2,...,q}\|\widehat{\mathcal{Q}}_{(k)}^{-\frac{1}{2}}\widehat{\mathcal{A}}^{H}_{(k)}S_{k_{i}}\|_{F}^{2}}\right)^{t}\|\mathcal{X}^{0}-\mathcal{X}^{\star}\|_{F(\mathcal{Q})}^{2}.
\end{align*}
\end{corollary}
\emph{Proof:}
 This proof is similar to that of Corollary \ref{sec3.2cor3.1}, so we omit it  here.

 The other improved strategy is to take the real part of the complex approximate solution directly. The specific algorithm is presented in Algorithm \ref{TSP-2}. It 
 has good performance confirmed by numerical experiments in Section \ref{sec:exp}. Unfortunately, we can't provide its rigorous theoretical analysis. 

\begin{algorithm}[htbp]
\caption{TSP-\uppercase\expandafter{\romannumeral2} method}\label{TSP-2}
$\textbf{Input:}$ $\mathcal{X}^{0}\in \mathbb{K}^{n\times p}_{l}$, $\mathcal{A}\in \mathbb{K}^{m\times n}_{l}$, $\mathcal{B}\in \mathbb{K}^{m\times p}_{l}$

$\textbf{Parameters:}$ fixed distribution $\mathfrak{D}_{k}$ over random matrices for $k=1,2,\cdots,l$, T-symmetric T-positive definite tubal matrix $\mathcal{Q}\in \mathbb{K}^{n\times n}_{l}$

$\widehat{\mathcal{X}}^{0}\longleftarrow \text{fft}(\mathcal{X}^{0},[~],3)$, $\widehat{\mathcal{A}}\longleftarrow \text{fft}(\mathcal{A},[~],3)$, $\widehat{\mathcal{B}}\longleftarrow \text{fft}(\mathcal{B},[~],3)$,$\widehat{\mathcal{Q}}\longleftarrow \text{fft}(\mathcal{Q},[~],3)$

$\textbf{for}$ $t=0,1,2,\cdots$

~~~~~$\textbf{for}$ $k=1,2,\cdots,l$

~~~~~~~~~$S_{k}\sim \mathfrak{D}_{k}$

~~~~~~~~~Compute $\widehat{\mathcal{G}}_{(k)}=S_{k}\left(S_{k}^{H}\widehat{\mathcal{A}}_{(k)}{\widehat{\mathcal{Q}}_{(k)}}^{-1}{\widehat{\mathcal{A}}_{(k)}}^{H}S_{k}\right)^{\dag}S_{k}^{H}$

~~~~~~~~~$\widehat{\mathcal{X}}^{t+1}_{(k)}=\widehat{\mathcal{X}}^{t}_{(k)}-{\widehat{\mathcal{Q}}_{(k)}}^{-1}{\widehat{\mathcal{A}}_{(k)}}^{H}\widehat{\mathcal{G}}_{(k)}\left(\widehat{\mathcal{A}}_{(k)}\widehat{\mathcal{X}}_{(k)}^{t}-\widehat{\mathcal{B}}_{(k)}\right)$

~~~~~$\textbf{end for}$

$\textbf{end for}$

$\mathcal{X}^{t+1}\longleftarrow \mathbf{Re}\left(\text{ifft}\left(\widehat{\mathcal{X}}^{t+1},[~],3\right)\right)$

 $\textbf{Output:}$ last iterate $\mathcal{X}^{t+1}$
\end{algorithm}

\begin{remark}\rm 
 Both the two improved strategies can be combined with the NTSP and three adaptive TSP methods. They can be called  NTSP-\uppercase\expandafter{\romannumeral1},
 ATSP-MD-\uppercase\expandafter{\romannumeral1}, ATSP-PR-\uppercase\expandafter{\romannumeral1}, ATSP-CS-\uppercase\expandafter{\romannumeral1},
 NTSP-\uppercase\expandafter{\romannumeral2}, ATSP-MD-\uppercase\expandafter{\romannumeral2}, ATSP-PR-\uppercase\expandafter{\romannumeral2} and ATSP-CS-\uppercase\expandafter{\romannumeral2} methods. The details of these algorithms are omitted here.
\end{remark}

\begin{remark}\rm
For $k=1,2,\cdots,l$, if we choose $S_{k}=e_{k_{i}}\in \mathbb{R}^{m}$ (the unit coordinate vector in $\mathbb{R}^{m}$) with $i=1,2,\cdots,m$ and $\mathcal{Q}=\mathcal{I} \in \mathbb{K}^{n\times n}_{l}$ in Algorithm \ref{TSP-1} and \ref{TSP-2}, we can obtain two improved TRK (i.e., TRK-\uppercase\expandafter{\romannumeral1} and TRK-\uppercase\expandafter{\romannumeral2}) methods. The convergence guarantees of the TRK-\uppercase\expandafter{\romannumeral1} method can be obtained according to Corollary \ref{sec3.3cor3.2}. That is, when selecting $k_{i}$ with probability proportional to the magnitude of row $k_{i}$ of $\widehat{\mathcal{A}}_{(k)}$, we have
\begin{align*}
\mathbf{E}\left[\left\|\mathcal{X}^{t}-\mathcal{X}^{\star}\right\|_{F}^{2}|\mathcal{X}^{0}\right]\leq\left(1-\mathop{\min}_{k=1,2,\cdots,l}\frac{\lambda_{\min}\left(\widehat{\mathcal{A}}_{(k)}\widehat{\mathcal{A}}_{(k)}^{H}\right)}{\|\widehat{\mathcal{A}}_{(k)}\|_{F}^{2}}\right)^{t}\|\mathcal{X}^{0}-\mathcal{X}^{\star}\|_{F}^{2}.
\end{align*}
\end{remark}

\section{Numerical experiments}\label{sec:exp}
\subsection{Implementation tricks and computation complexity}
In this subsection, we discuss the computation costs at each iteration of some nonadaptive and adaptive TSP methods proposed in previous sections. Specifically, the nonadaptive methods include the NTSP, NTSP-\uppercase\expandafter{\romannumeral1} and NTSP-\uppercase\expandafter{\romannumeral2} methods, and the adaptive methods include the ATSP-MD, ATSP-PR, ATSP-CS,  ATSP-MD-\uppercase\expandafter{\romannumeral2}, ATSP-PR-\uppercase\expandafter{\romannumeral2} and ATSP-CS-\uppercase\expandafter{\romannumeral2} methods. Similar to \cite{gower2019adaptive}, we implement these methods except the NTSP-\uppercase\expandafter{\romannumeral1} one in their corresponding fast versions in the following numerical experiments, for example, Algorithm \ref{EATSP} and \ref{EAITSP} are the fast versions of the  ATSP-PR and  ATSP-PR-\uppercase\expandafter{\romannumeral2} methods, respectively.
\begin{algorithm}[htbp]
\caption{Fast ATSP-PR method in Fourier domain}\label{EATSP}
$\textbf{Input:}$ $\mathcal{X}^{0}\in \mathbb{K}^{n\times p}_{l}$, $\mathcal{A}\in \mathbb{K}^{m\times n }_{l}$, and $\mathcal{B}\in \mathbb{K}^{m\times p }_{l}$

$\textbf{Parameters:}$ a set of sketching tubal matrices $\boldsymbol{\mathcal{S}}=[\mathcal{S}_{1},\cdots,\mathcal{S}_{q}]$, T-symmetric T-positive definite tubal matrix $\mathcal{Q}\in \mathbb{K}^{n\times n}_{l}$

\textbf{1}: $\widehat{\mathcal{X}}^{0}\longleftarrow \text{fft}(\mathcal{X}^{0},[~],3)$, $\widehat{\mathcal{A}}\longleftarrow
\text{fft}(\mathcal{A},[~],3)$, $\widehat{\mathcal{B}}\longleftarrow \text{fft}(\mathcal{B},[~],3)$, $\widehat{\mathcal{Q}}\longleftarrow \text{fft}(\mathcal{Q},[~],3)$, $\widehat{\mathcal{S}_{i}}\longleftarrow \text{fft}(\mathcal{S}_{i},[~],3)$ for $i=1,2,\cdots,q$

\textbf{2}: $\textbf{for}$ $k=1,2,\cdots,l$

\textbf{3}: ~~~~~Compute $(\widehat{\mathcal{C}}_{i})_{(k)}=\text{Cholesky}\left(\left((\widehat{\mathcal{S}}_{i})_{(k)}^{H}\widehat{\mathcal{A}}_{(k)}\widehat{\mathcal{Q}}^{-1}_{(k)}\widehat{\mathcal{A}}_{(k)}^{H}(\widehat{\mathcal{S}_{i}})_{(k)}\right)^{\dag}\right)$ for $i=1,2,\cdots,q$

\textbf{4}: ~~~~~Compute $\widehat{\mathcal{Q}}^{-1}_{(k)}\widehat{\mathcal{A}}_{(k)}^{H}(\widehat{\mathcal{S}_{i}})_{(k)}(\widehat{\mathcal{C}}_{i})_{(k)}$ for $i=1,2,\cdots,q$

\textbf{5}: ~~~~~Compute $(\widehat{\mathcal{C}}_{i})_{(k)}^{H}(\widehat{\mathcal{S}_{i}})_{(k)}^{H}\widehat{\mathcal{A}}_{(k)}\widehat{\mathcal{Q}}^{-1}_{(k)}\widehat{\mathcal{A}}_{(k)}^{H}(\widehat{\mathcal{S}}_{j})_{(k)}(\widehat{\mathcal{C}}_{j})_{(k)}$ for $i,j=1,2,\cdots,q$

\textbf{6}: ~~~~~Initialize $(\widehat{\mathcal{R}}_{i}^{0})_{(k)}=(\widehat{\mathcal{C}}_{i})_{(k)}^{H}\left((\widehat{\mathcal{S}_{i}})_{(k)}^{H}\left(\widehat{\mathcal{A}}_{(k)}\widehat{\mathcal{X}}_{(k)}^{0}-\widehat{\mathcal{B}}_{(k)}\right)\right)$ for $i=1,2,\cdots,q$

\textbf{7}: $\textbf{end for}$

\textbf{8}: $\textbf{for}$ $t=0,1,2,\cdots$

\textbf{9}:  ~~~~~~$f_{i}(\mathcal{X}^{t})=(1/l)\sum_{k=1}^{l}\|(\widehat{\mathcal{R}}_{i}^{t})_{(k)}\|_{F}^{2}$ for $i=1,2,\cdots,q$

\textbf{10}: ~~~~~Sample $i^{t} \sim \mathbf{p}^{t}$, where $p_{i}^{t}=f_{i}(\mathcal{X}^{t})/(\sum_{i=1}^{q}f_{i}(\mathcal{X}^{t}))$ for $i=1,2,\cdots,q$

\textbf{11}: ~~~~~$\textbf{for}$ $k=1,2,\cdots,l$

\textbf{12}: ~~~~~~~~~~Update $\widehat{\mathcal{X}}^{t+1}_{(k)}=\widehat{\mathcal{X}}^{t}_{(k)}-\left(\widehat{\mathcal{Q}}^{-1}_{(k)}\widehat{\mathcal{A}}_{(k)}^{H}(\widehat{\mathcal{S}}_{i^{t}})_{(k)}(\widehat{\mathcal{C}}_{i^{t}})_{(k)}\right)(\widehat{\mathcal{R}}_{i^{t}}^{t})_{(k)}$

\textbf{13}: ~~~~~~~Update $(\widehat{\mathcal{R}}_{i}^{t+1})_{(k)}=(\widehat{\mathcal{R}}_{i}^{t})_{(k)}-\left((\widehat{\mathcal{C}}_{i})_{(k)}^{H}(\widehat{\mathcal{S}_{i}})_{(k)}^{H}\widehat{\mathcal{A}}_{(k)}\widehat{\mathcal{Q}}^{-1}_{(k)}\widehat{\mathcal{A}}_{(k)}^{H}(\widehat{\mathcal{S}}_{i^{t}})_{(k)}(\widehat{\mathcal{C}}_{i^{t}})_{(k)}\right)(\widehat{\mathcal{R}}_{i^{t}}^{t})_{(k)}$ for $i=1,2,\cdots,q$

\textbf{14}: ~~~~~$\textbf{end for}$

\textbf{15}: $\textbf{end for}$

\textbf{16}: $\mathcal{X}^{t+1}\longleftarrow \text{ifft}\left(\widehat{\mathcal{X}}^{t+1},[~],3\right)$

$\textbf{Output:}$ last iterate $\mathcal{X}^{t+1}$

\end{algorithm}

\begin{algorithm}[htbp]
\caption{Fast ATSP-PR-\uppercase\expandafter{\romannumeral2} method in Fourier domain}\label{EAITSP}
$\textbf{Input:}$ $\mathcal{X}^{0}\in \mathbb{K}^{n\times p}_{l}$, $\mathcal{A}\in \mathbb{K}^{m\times n }_{l}$, and $\mathcal{B}\in \mathbb{K}^{m\times p }_{l}$

$\textbf{Parameters:}$ $k$ sets of sketching matrices $\boldsymbol{S}_{k}=[S_{k_{1}},\cdots,\mathcal{S}_{k_{q}}]$ for $k=1,2,\cdots,l$, T-symmetric T-positive definite tubal matrix $\mathcal{Q}\in \mathbb{K}^{n\times n}_{l}$

\textbf{1}: $\widehat{\mathcal{X}}^{0}\longleftarrow \text{fft}(\mathcal{X}^{0},[~],3)$, $\widehat{\mathcal{A}}\longleftarrow
\text{fft}(\mathcal{A},[~],3)$, $\widehat{\mathcal{B}}\longleftarrow \text{fft}(\mathcal{B},[~],3)$, $\widehat{\mathcal{Q}}\longleftarrow \text{fft}(\mathcal{Q},[~],3)$

\textbf{2}: $\textbf{for}$ $k=1,2,\cdots,l$

\textbf{3}: ~~~~~Compute $C_{k_{i}}=\text{Cholesky}\left(\left(S_{k_{i}}^{H}\widehat{\mathcal{A}}_{(k)}\widehat{\mathcal{Q}}^{-1}_{(k)}\widehat{\mathcal{A}}_{(k)}^{H}S_{k_{i}}\right)^{\dag}\right)$ for $i=1,2,\cdots,q$

\textbf{4}: ~~~~~Compute $\widehat{\mathcal{Q}}^{-1}_{(k)}\widehat{\mathcal{A}}_{(k)}^{H}S_{k_{i}}C_{k_{i}}$ for $i=1,2,\cdots,q$

\textbf{5}: ~~~~~Compute $C_{k_{i}}^{H}S_{k_{i}}^{H}\widehat{\mathcal{A}}_{(k)}\widehat{\mathcal{Q}}^{-1}_{(k)}\widehat{\mathcal{A}}_{(k)}^{H}S_{k_{j}}C_{k_{j}}$ for $i,j=1,2,\cdots,q$

\textbf{6}: ~~~~~Initialize $R_{k_{i}}^{0}=C_{k_{i}}^{H}\left(S_{k_{i}}^{H}\left(\widehat{\mathcal{A}}_{(k)}\widehat{\mathcal{X}}_{(k)}^{0}-\widehat{\mathcal{B}}_{(k)}\right)\right)$ for $i=1,2,\cdots,q$

\textbf{7}: $\textbf{end for}$

\textbf{8}: $\textbf{for}$ $t=0,1,2,\cdots$

\textbf{9}: ~~~~~$\textbf{for}$ $k=1,2,\cdots,l$

\textbf{10}: ~~~~~~~~~~$f_{i}(\widehat{\mathcal{X}}_{(k)}^{t})=\|R_{k_{i}}^{t}\|_{F}^{2}$ for $i=1,2,\cdots,q$

\textbf{11}: ~~~~~~~~~~Sample $k_{i}^{t} \sim\mathbf{p}_{k}^{t}$, where $p_{k_{i}}^{t}=f_{i}(\widehat{\mathcal{X}}_{(k)}^{t})/(\sum_{i=1}^{q}f_{i}(\widehat{\mathcal{X}}_{(k)}^{t})$ for $i=1,2,\cdots,q$

\textbf{12}: ~~~~~~~~~~Update $\widehat{\mathcal{X}}^{t+1}_{(k)}=\widehat{\mathcal{X}}^{t}_{(k)}-\left(\widehat{\mathcal{Q}}^{-1}_{(k)}\widehat{\mathcal{A}}_{(k)}^{H}S_{k_{i}^{t}}C_{k_{i}^{t}}\right)R_{k_{i}^{t}}^{t}$

\textbf{13}: ~~~~~~~~~~Update $R_{k_{i}}^{t+1}=R_{k_{i}}^{t}-\left(C_{k_{i}}^{H}S_{k_{i}}^{H}\widehat{\mathcal{A}}_{(k)}\widehat{\mathcal{Q}}^{-1}_{(k)}\widehat{\mathcal{A}}_{(k)}^{H}S_{k_{i}^{t}}C_{k_{i}^{t}}\right)R_{k_{i}^{t}}^{t}$ for $i=1,2,\cdots,q$

\textbf{14}: ~~~~~$\textbf{end for}$

\textbf{15}: $\textbf{end for}$

\textbf{16}: $\mathcal{X}^{t+1}\longleftarrow \mathbf{Re}\left(\text{ifft}\left(\widehat{\mathcal{X}}^{t+1},[~],3\right)\right)$

$\textbf{Output:}$ last iterate $\mathcal{X}^{t+1}$
\end{algorithm}

We first consider the computation complexities of the 
NTSP, ATSP-MD, ATSP-PR and ATSP-CS methods. Since the difference of the fast versions of these methods mainly lies in 
how to compute the sampling probabilities, we first present the flops of each step of the four algorithms without the step on sampling:

1. Computing the sketched losses $\{f_{i}(\mathcal{X}^{t}): i=1,2,\cdots,q\}$ requires $2\tau plq~(l>1)$ or $(2\tau p-1)q~(l=1)$ flops if the sketched residuals $\{\widehat{\mathcal{R}}^{t}_{i}: i=1,2,\cdots,q\}$ are precomputed. 

2. Updating $\widehat{\mathcal{X}}^{t}$ to $\widehat{\mathcal{X}}^{t+1}$ requires $2\tau npl$ flops when $$\{\widehat{\mathcal{Q}}^{-1}_{(k)}\widehat{\mathcal{A}}_{(k)}^{H}(\widehat{\mathcal{S}_{i}})_{(k)}(\widehat{\mathcal{C}}_{i})_{(k)}: i=1,2,\cdots,q, k=1,2,\cdots,l\}$$ are precomputed. 

3. Updating $\{\widehat{\mathcal{R}}^{t}_{i}: i=1,2,\cdots,q\}$ to $\{\widehat{\mathcal{R}}^{t+1}_{i}: i=1,2,\cdots,q\}$ requires $2\tau^{2}plq$ flops if $$\{(\widehat{\mathcal{C}}_{i})_{(k)}^{H}(\widehat{\mathcal{S}_{i}})_{(k)}^{H}\widehat{\mathcal{A}}_{(k)}\widehat{\mathcal{Q}}^{-1}_{(k)}\widehat{\mathcal{A}}_{(k)}^{H}(\widehat{\mathcal{S}}_{j})_{(k)}(\widehat{\mathcal{C}}_{j})_{(k)} : i,j=1,2,\cdots,q, k=1,2,\cdots,l\}$$ are precomputed. Note that for the NTSP method, one only needs to compute the single sketched residual $\widehat{\mathcal{R}}_{i^{t}}^{t}$, where $(\widehat{\mathcal{R}}_{i^{t}}^{t})_{(k)}=(\widehat{\mathcal{C}}_{i^{t}})_{(k)}^{H}\left((\widehat{\mathcal{S}}_{i^{t}})_{(k)}^{H}\left(\widehat{\mathcal{A}}_{(k)}\widehat{\mathcal{X}}_{(k)}^{t}-\widehat{\mathcal{B}}_{(k)}\right)\right)$ for $k=1,2,\cdots,l$. If $(\widehat{\mathcal{C}}_{i})_{(k)}^{H}(\widehat{\mathcal{S}_{i}})_{(k)}^{H}\widehat{\mathcal{A}}_{(k)}$ and $(\widehat{\mathcal{C}}_{i})_{(k)}^{H}(\widehat{\mathcal{S}_{i}})_{(k)}^{H}\widehat{\mathcal{B}}_{(k)}$ are precomputed for $i=1,2,\cdots,q$, $k=1,2,\cdots,l$, computing sketched residual $\widehat{\mathcal{R}}_{i^{t}}^{t}$ directly from the iterate $\mathcal{X}^{t}$ costs $2\tau npl$ flops. Hence, when $\tau q>n$, it is cheaper for the NTSP method to compute the sketched residual $\widehat{\mathcal{R}}_{i^{t}}^{t}$ directly than using update  formula.

Therefore, the nonsampling flops of the NTSP method and the adaptive cases (ATSP-MD, ATSP-PR, ATSP-CS) are $2\tau pl \min(n,\tau q)+2\tau npl$ and $(2\tau^{2} p+2\tau p)lq+2\tau npl $ ($l>1$) or $(2\tau^{2} p+2\tau p -1)q+2\tau np $ ($l=1$), respectively. 

Next, we give the cost of computing the sampling probabilities $\mathbf{p}^{t}$ from the sketched losses $\{f_{i}(\mathcal{X}^{t}): i=1,2,\cdots,q\}$. It depends on the sampling strategy. Specifically, 
for the NTSP method, it requires $\mathcal{O}(1)$ flops; for the ATSP-MD method, it needs $q$ flops if $\tau>1$ and $\mathcal{O}(log(q))$ flops if $\tau =1$; the ATSP-PR method requires approximately $2q$ flops on average; the ATSP-CS method requires $6q$ flops. 

Putting all the costs together, we report the total costs per iteration of the above four methods 
in Table \ref{tab:my_label3}.

\begin{table}
    \caption{The computation costs for the nonadaptive and adaptive TSP methods, where $\tau$ is the sketch size, $q$ is the number of sketches, $n$ and $l$ are the dimension of $\mathcal{A}$, and $p$ is the size of $\mathcal{B}$.}
    \centering
    \begin{tabular}{|c|c|c|}
        \hline
         Method & Flops per iteration when $\tau >1$ & Flops per iteration when $\tau =1$  \\
        \hline
         NTSP  & $2\tau pl \min(n,\tau q)+2\tau npl$ & $2pl\min(n,q)+2npl$   \\
        \hline
        \multirow{2}{*}{ATSP-MD} & $(2\tau^{2} pl+2\tau pl +1 )q+2\tau npl$ if $l>1$ &  $ 4plq+2npl$ if $l>1$ \\ 
         & $(2\tau^{2}p+2\tau p)q+2\tau np$ if $l=1$ & $(4p-1)q+2np$ if $l=1$\\
        \hline
        \multirow{2}{*}{ATSP-PR} &  $(2\tau^{2} pl+2\tau pl +2)q+2\tau npl$ if $l>1$ &  $(4pl+2)q+2npl$ if $l>1$ \\ & $(2\tau^{2}p+2\tau p+1)q+2\tau np$ if $l=1$ &  $(4p+1)q+2np$ if $l=1$ \\
        \hline
         \multirow{2}{*}{ATSP-CS} &  $(2\tau^{2} pl+2\tau pl +6)q+2\tau npl$ if $l>1$ &  $(4pl+6)q+2npl$ if $l>1$ \\ & $(2\tau^{2}p+2\tau p+5)q+2\tau np$ if $l=1$ &  $(4p+5)q+2np$ if $l=1$ \\
        \hline
    \end{tabular}
    
    \label{tab:my_label3}
\end{table}

In a similar way, we can give the computation costs at each iteration of the NTSP-\uppercase\expandafter{\romannumeral2}, ATSP-MD-\uppercase\expandafter{\romannumeral2}, ATSP-PR-\uppercase\expandafter{\romannumeral2} and ATSP-CS-\uppercase\expandafter{\romannumeral2} methods. The details are omitted here, and the total costs per iteration 
are reported in Table \ref{tab:my_label6}.

\begin{table}
    \caption{The computation costs of the improved nonadaptive and adaptive TSP methods, where $\tau$ is the sketch size, $q$ is the number of sketches, $n$ and $l$ are the dimension of $\mathcal{A}$, and $p$ is the size of $\mathcal{B}$.}
    \centering
    \begin{tabular}{|c|c|c|}
        \hline
         Method & Flops per iteration when $\tau >1$ & Flops per iteration when $\tau =1$  \\
        \hline
         NTSP-\uppercase\expandafter{\romannumeral2}    & $\mathcal{O}(\tau pln)$ & $\mathcal{O}(pln)$   \\
        \hline
        ATSP-MD-\uppercase\expandafter{\romannumeral2}  & $(2\tau^{2} p+2\tau p )ql+2\tau npl$ & $\mathcal{O}(\max(q,n)pl)$ \\
        \hline
        ATSP-PR-\uppercase\expandafter{\romannumeral2}  &$(2\tau^{2} p+2\tau p +1)ql+2\tau npl$ & $(4p+1)ql+2npl$ \\
        \hline
        ATSP-CS-\uppercase\expandafter{\romannumeral2}  &$(2\tau^{2} p+2\tau p +5)ql+2\tau npl$ & $(4p+5)ql+2npl$ \\
        \hline
    \end{tabular}
    
    \label{tab:my_label6}
\end{table}

For the NTSP-\uppercase\expandafter{\romannumeral1} method, it has no fast implement version. We present separately its complexities of each step as follows:

1. Computing $\widehat{\mathcal{A}_{\mathcal{S}}}$ and $\widehat{\mathcal{B}_{\mathcal{S}}}$ requires $\mathcal{O}(\tau nml+\tau nl\log l)$ and $\mathcal{O}(\tau pml+\tau pl\log l)$ flops, respectively. If \{$S_{k}\in\mathbb{C}^{m\times \tau}$: $k=1,2,...,l$\} are random sampling matrices, then copmputing $\widehat{\mathcal{A}_{\mathcal{S}}}$ and $\widehat{\mathcal{B}_{\mathcal{S}}}$ requires $\mathcal{O}(\tau nl\log l)$ and $\mathcal{O}(\tau pl\log l)$ flops, respectively.

2. Updating $\widehat{\mathcal{X}}^{t}$ to $\widehat{\mathcal{X}}^{t+1}$ requires $\mathcal{O}(n^{2}\tau l +np\tau l+n\tau^{2}l+\tau^{3}l)$ flops when \{$\widehat{\mathcal{Q}}_{(k)}^{-1}: k=1,2,...,l$\} is precomputed. Note that, if $\mathcal{Q}=\mathcal{I} \in \mathbb{K}^{n\times n}_{l}$ , then updating $\widehat{\mathcal{X}}^{t}$ to $\widehat{\mathcal{X}}^{t+1}$ requires $\mathcal{O}(\tau npl+n\tau^{2}l+\tau^{3}l)$ flops.

\subsection{Examples}
We use four numerical experiments to illustrate the performance of the proposed TSP method and its adaptive variants for solving the tensor linear systems (\ref{sec1e1}). To compare with the existing methods more intuitively, we only consider the relevant experiments on a special case of the TSP method, i.e., the TRK method. Specifically, we compare the performance of ten algorithms including four nonadaptive TRK methods, i.e., NTRKU (uniform sampling),\cite{ma2021randomized} NTRKS (probabilities proportional to the magnitude of horizontal slices of $\mathcal{A}$),\cite{ma2021randomized} NTRKR-\uppercase\expandafter{\romannumeral1} (probabilities proportional to the magnitude of the rows of the frontal slices of $\widehat{\mathcal{A}}$) and NTRKR-\uppercase\expandafter{\romannumeral2} (probabilities proportional to the magnitude of the rows of the frontal slices of $\widehat{\mathcal{A}}$), as well as six adaptive methods, i.e., ATRKS-PR, ATRKS-MD, ATRKS-CS, ATRKR-PR-\uppercase\expandafter{\romannumeral2}, ATRKR-MD-\uppercase\expandafter{\romannumeral2} and ATRKR-CS-\uppercase\expandafter{\romannumeral2}. The relative error used to determine the effectiveness of these different methods is defined as $$\varepsilon=\frac{\|\mathcal{X}^{t}-\mathcal{X}^{\star}\|_{F}}{\|\mathcal{X}^{\star}\|_{F}}.$$ 
We run each method until the relative error is below $10^{-10}$ (Example \ref{E1}), $10^{-4}$ (Example \ref{E2} and \ref{E3}) or $0.005$ (Example \ref{E4}). In the following examples, we use $\mathcal{X}^{0}=\emph{O}$ as an initial point and all results are average on 10 trails. In each figure, we plot
the relative error (i.e., Error) on the vertical axis, starting with 1. For the horizontal axis, we use either the number of the iterations (i.e., Iters) or running
time (i.e., Time(s)). 
Note that we do not consider the precomputational cost, but only the costs spent at each iteration. All computations were carried out in MATLAB R2018a on a standard MacBook Pro 2019 with an Intel Core i9 processor and 16GB memory.     

\begin{example}[synthetic data]\label{E1}\rm
Let the entries of $\mathcal{A}\in \mathbb{K}^{m\times n}_{l}$ and $\mathcal{X}\in \mathbb{K}^{n\times p}_{l}$ be drawn i.i.d. from a standard Gaussian distribution, and the right-hand tubal matrix be $\mathcal{B}=\mathcal{A}*\mathcal{X}\in \mathbb{K}^{m\times p}_{l}$. Specifically, we compare the empirical performance of the ten algorithms for a system with $m=500$, $n=200$, $p=50$ and $l=50$. Figure \ref{fig1} 
shows that the NTRKR-\uppercase\expandafter{\romannumeral1} method has the best performance in terms of CPU time among the ten methods, and has the fewest iteration steps among the four nonadaptive methods. The other three nonadaptive methods (i.e., NTRKU, NTRKS and NTRKR-\uppercase\expandafter{\romannumeral2}) perform similarly. The number of iteration steps of each of the six adaptive methods is smaller than that of each of the four nonadaptive methods, and, except for NTRKR-\uppercase\expandafter{\romannumeral1}, the time of each of the six adaptive methods is also less than that of each of the other three nonadaptive methods. In addition, compared with the original adaptive methods (i.e., ATRKS-PR, ATRKS-MD and ATRKS-CS), the adaptive methods combined with the second improved strategy (i.e., ATRKR-PR-\uppercase\expandafter{\romannumeral2}, ATRKR-MD-\uppercase\expandafter{\romannumeral2} and ATRKR-CS-\uppercase\expandafter{\romannumeral2}) vastly reduced the number of iteration steps and CPU time. 

\begin{figure}[tbhp] 
\centering
\includegraphics[width=1 \textwidth]{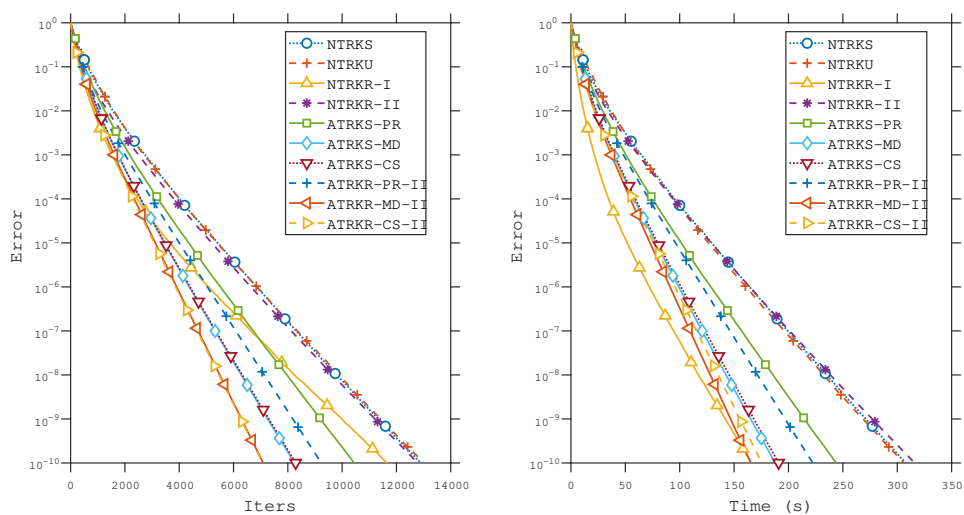}
\caption{ Errors versus iterations (left) and CPU time (right) for four nonadaptive 
and six adaptive \texttt{TRK} methods on synthetic data with $m=500$, $n=200$, $p=50$ and $l=50$.}\label{fig1}
\end{figure}

\end{example}

\begin{example}[CT data]\label{E2}\rm
In this experiment, we evaluate the performance of the ten methods on real world CT data set. The underlying signal $\mathcal{X}$ is a tubal matrix of size $512\times 512 \times 11$, where each frontal slice is a $512\times 512$ matrix of the C1-vertebrae. The images for the experiment were obtained from the Laboratory of the Human Anatomy and Embryology, University of Brussels (ULB), Belgium.\cite{ctdata} To set up the tensor linear system, we generate randomly a Gaussian tubal matrix $\mathcal{A}\in \mathbb{K}^{1000\times 512}_{11}$ and form the measurement tubal matrix $\mathcal{B}$ by $\mathcal{B}=\mathcal{A}*\mathcal{X}$. The numerical results of this experiment are provided in Figure \ref{fig2}, from which we can see that the performance of the four nonadaptive methods is almost the same except that the NTRKR-\uppercase\expandafter{\romannumeral1} method takes less time. 
Among the ten methods, the six adaptive methods outperform the four nonadaptive methods in terms of iteration numbers. For running time, they perform better than the NTRKU, NTRKS, and NTRKR-II methods. 
While, the NTRKR-\uppercase\expandafter{\romannumeral1} method spends less running time than two adaptive methods, i.e., ATRKS-PR and ATRKR-PR-\uppercase\expandafter{\romannumeral2}. In addition, like Example \ref{E1}, the adaptive methods combined with the second improved strategy perform better than the original adaptive methods. 

\begin{figure}[tbhp] 
\centering
\includegraphics[width=1 \textwidth]{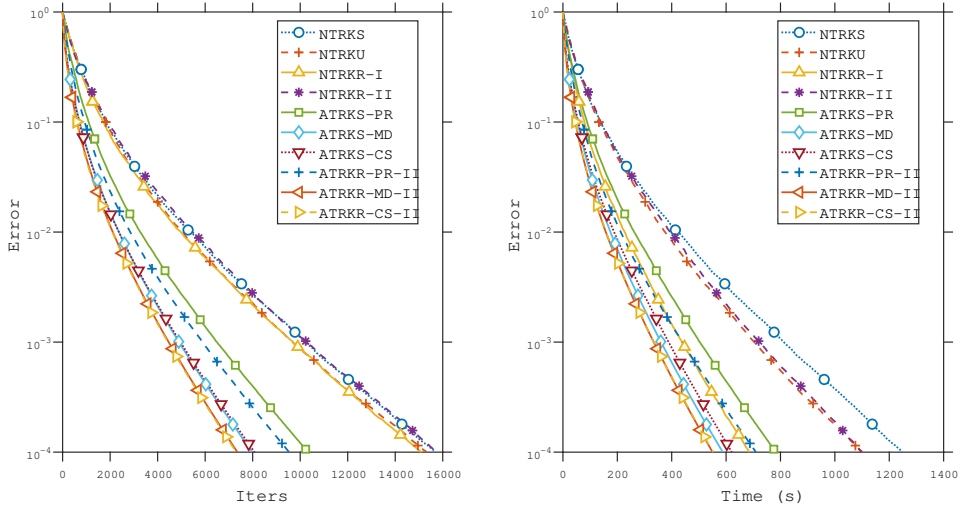}
\caption{Errors versus iterations (left) and CPU time (right) for four nonadaptive 
and six adaptive \texttt{TRK} methods on real world CT data set.}\label{fig2}
\end{figure}

\end{example}

\begin{example}[video data]\label{E3}\rm
The example illustrates the performance of the ten methods on the video data where the frontal slices of the tubal matrix $\mathcal{X}$ are the first $80$ frames from the 1929 film "Finding His Voice".\cite{videodata} Each video frame has $480\times 368$ pixels. Similar to Example \ref{E2}, we generate randomly a Gaussian tubal matrix $\mathcal{A}\in \mathbb{K}^{1000\times 480}_{80}$ and form the measurement tubal matrix $\mathcal{B}$ by $\mathcal{B}=\mathcal{A}*\mathcal{X}$. From Figure \ref{fig3}, we can find that for the four nonadaptive methods, they are similar in the number of iteration steps, however, in terms of CPU time, the NTRKR-\uppercase\expandafter{\romannumeral1} method is the fastest one and the NTRKR-\uppercase\expandafter{\romannumeral2} method is faster than the NTRKS and NTRKU methods. Among the ten methods, the six adaptive methods outperform the four nonadaptive methods in terms of iteration 
numbers. For running time, they perform better than the NTRKU, NTRKS, and NTRKR-\uppercase\expandafter{\romannumeral2} methods. 
While, only two adaptive methods, i.e., ATRKS-CS and ATRKR-CS-\uppercase\expandafter{\romannumeral2}, are faster than the NTRKR-\uppercase\expandafter{\romannumeral1} method. In addition, as in the previous examples, combining with the second improved strategy can indeed improve the adaptive methods in terms of the number of iteration steps and CPU time. 

\begin{figure}[tbhp] 
\centering
\includegraphics[width=1 \textwidth]{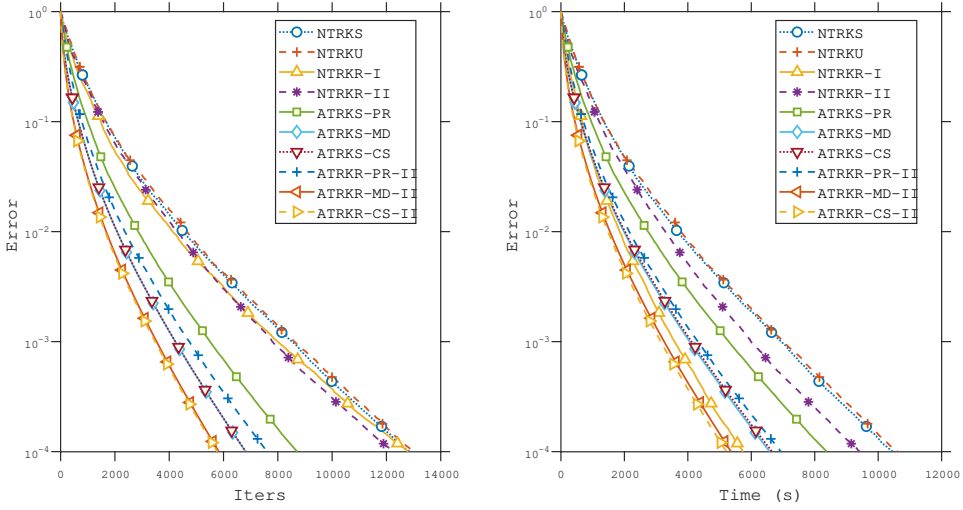}
\caption{Errors versus iterations (left) and CPU time (right) for four nonadaptive 
and six adaptive \texttt{TRK} methods on video data set.}\label{fig3}
\end{figure}

\end{example}

\begin{example}[image deblurring]\label{E4}\rm
This example considers an image sequence $\{\widehat{X}_{j}\}_{j=1}^{27}$ from a 3D MRI image data set \textbf{mri} in MATLAB, which has $27$ slices with dimensions $128\times 128$. Assume that each image is degraded by a Gaussian convolution kernel $\widehat{H}$ of size $5\times 5$ with standard deviation $2$.  By the construction, we can obtain the image deblurring problem as follows:
\begin{align}\label{sec6.e1}
    H \circledast X_{j}=Y_{j},~~\text{for}~~j=1,2,\cdots, 27,
\end{align}
where $\{X_{j}\}_{j=1}^{27}$ and $H$ are extended by padding $\{\widehat{X}_{j}\}_{j=1}^{27}$ and $\widehat{H}$ with the zeros respectively, and they are of size $132\times 132$; $Y_{j}\in \mathbb{R}^{132\times 132}$, for $j=1,2,\cdots,27$, are the observed blurry images; and $\circledast$ is the 2D convolution. Using the equivalence between 2D convolution and t-product, the above problem (\ref{sec6.e1}) can be equivalently rewritten as the following tensor linear system
\begin{align}\label{sec6.e2}
 \mathcal{A}*\mathcal{X}=\mathcal{B},
\end{align}
where $\mathcal{A}\in \mathbb{K}^{132\times 132}_{132}$ whose $k$-th frontal slice is the circulant matrix generated by the $k$-th column of $H$, i.e., $\mathcal{A}_{(k)}=\text{circ}\left(H_{(:,k)}\right)$ for $k=1,2,\cdots,132$; $\mathcal{X}$, $\mathcal{B}\in \mathbb{K}^{132\times 27}_{132}$ are the tubal matrices by setting $\mathcal{X}_{(i,j,k)}=(X_{j})_{(k,i)}$ and $\mathcal{B}_{(i,j,k)}=(Y_{j})_{(k,i)}$ for $i=1.2,\cdots,132$, $j=1,2,\cdots,27$ and $k=1,2,\cdots,132$, respectively. As shown in Figure \ref{fig4}, for the four nonadaptive methods, the NTRKU, NTRKS, and NTRKR-\uppercase\expandafter{\romannumeral2} methods have similar numerical performance, while the NTRKR-\uppercase\expandafter{\romannumeral1} method is better than the previous three methods in terms of the number of iteration steps and computing time. Except for the NTRKR-\uppercase\expandafter{\romannumeral1} method which is competitive with the six adaptive methods, the other three nonadaptive methods have considerably larger iteration numbers and more CPU time than the six adaptive methods. For the six adaptive methods, from the enlarged small graph in Figure \ref{fig4}, it can be found that the experimental performance is consistent with the previous numerical examples, that is, the combination with the second improved strategy can indeed make the adaptive methods perform better. In addition, the first slice of the clean image sequence, its corresponding blurry observation and the images recovered from the four nonadaptive and six adaptive methods are shown in Figure \ref{fig5}. 

\begin{figure}[tbhp] 
\centering
\includegraphics[width=1 \textwidth]{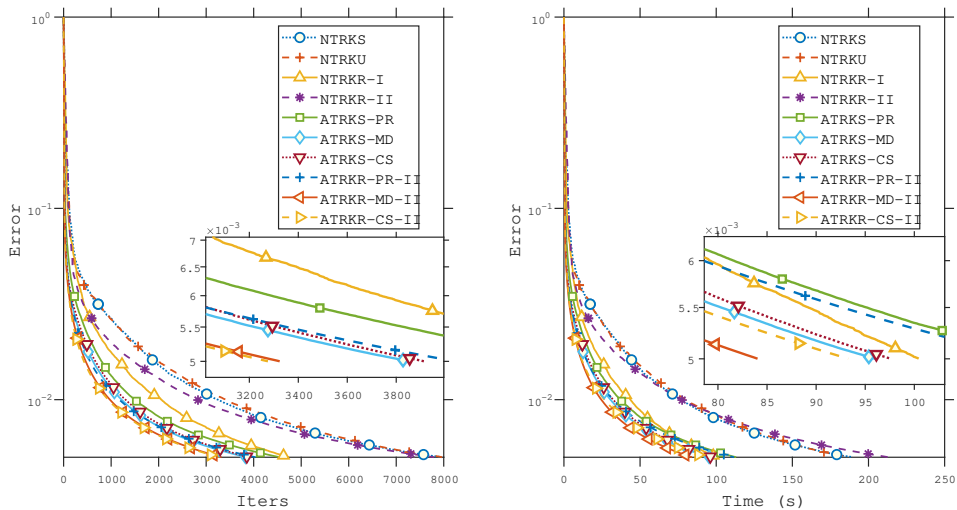}
\caption{Errors versus iterations (left) and CPU time (right) for four nonadaptive 
and six adaptive \texttt{TRK} methods on image deblurring problem.}\label{fig4}
\end{figure}

\begin{figure}[htbp]
\centering
\includegraphics[width=1 \textwidth]{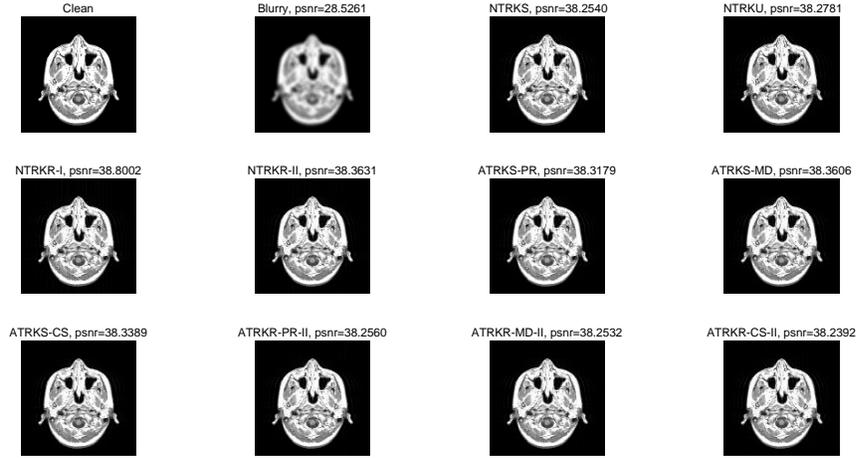}
\caption{The first slice of the clean image sequence, its corresponding blurry observation and the images recovered from the four nonadaptive 
and six adaptive \texttt{TRK} methods.}\label{fig5}
\end{figure}

\end{example}

\section{Conclusion}\label{sec:con}
In this paper, we propose the TSP method and its adaptive variants for tensor linear systems. We also discuss their Fourier domain versions. Two strategies used to improve the sketching or sampling techniques for real tensor linear systems are provided. Extensive numerical results including the ones from the CT signal recovery and image
deblurring problems show that the adaptive methods can indeed accelerate the nonadaptive ones and the two improved strategies are indeed effective for real linear systems. 

\bibliography{mybibfile}

\end{document}